\documentclass{amsart}
\usepackage{amsmath,amsfonts,amssymb,epsfig,graphicx,subfigure,graphs,amsopn,verbatim, color}

\newtheorem{theo}{Theorem}[section]
\newtheorem{prop}[theo]{Proposition}
\newtheorem{lemma}[theo]{Lemma}
\newtheorem{cor}[theo]{Corollary}
\newtheorem{defs}[theo]{Definition}
\numberwithin{equation}{section}
\theoremstyle{remark}
\newtheorem{rem}{Remark}

\newcommand{\Z}{\mathbb{Z}}

\newcommand{\N}{\mathbb{N}}
\newcommand{\C}{\mathbb{C}}
\newcommand{\B}{\mathcal{B}}
\newcommand{\ds}{\displaystyle}

\newcommand{\VE}{(\mathcal{V,E})}
\newcommand{\OVE}{(\mathcal{V,E},\geq)}
\newcommand{\vep}{\VE_{p(x)}}
\newcommand{\V}{\mathcal{V}}
\newcommand{\E}{\mathcal{E}}
\newcommand{\p}{\mathcal{P}}

\newcommand{\Sl}{\mathcal{S_L}}
\newcommand{\D}{\mathcal{D_L}}
\newcommand{\Slp}{(\Sl)_{p(x)}}
\newcommand{\Dlp}{(\D)_{p(x)}}

\newcommand{\X}{X_{p(x)}}
\newcommand{\T}{T_{p(x)}}
\newcommand{\co}{\text{coeff}_{p(x)}}
\newcommand{\G}{G_{p(x)}}

\begin{document}

\title{Limited scope adic transformations}

\author{Sarah Bailey Frick\\The Ohio State University, Department of Mathematics, 100 Math Tower, 231 W. 18th Ave., Columbus, OH 43210\\\email{frick@math.ohio-state.edu}}

\begin{abstract}
We introduce a family of adic transformations on diagrams that are nonstationary and nonsimple.  This family includes some previously studied adic transformations.  We relate the dimension group of each these diagrams to the dynamical system determined by the adic transformation on the infinite edge paths, and we explicitly compute the dimension group for a subfamily.  We also determine the ergodic adic invariant probability measures for this subfamily, and show that each system of the subfamily is loosely Bernoulli.  We also give examples of particular adic transformations with roots of unity as well as one which is totally ergodic called the Euler adic.  We also show that the Euler adic is loosely Bernoulli.  
\end{abstract}

\maketitle
\section{Introduction}\label{description}
  In 1972, Bratteli introduced infinite directed graphs known now as \emph{Bratteli diagrams} as a tool to study approximately finite-dimensional (AF) algebras, \cite{Bratteli}.  To each of these graphs, Vershik introduced a transformation, now known as the \emph{Bratteli-Vershik} or \emph{adic} transformation as a method of modeling cutting and stacking transformations, \cite{Vershik1981, Vershik1985, VerLiv1992}.  He also showed that every ergodic measure-preserving transformation on a Lebesgue space is isomorphic to a Bratteli-Vershik transformation which has a unique invariant measure associated to it.  Herman, Putnam, and Skau went on to show that every minimal homeomorphism of the Cantor set is topologically conjugate to a adic transformation with certain properites, \cite{HPS}.
  
  One adic transformation proposed by Vershik is the Pascal adic transformation.  This single transformation has been the subject of much study; see \cite{HIK, Vershik4, santiago, Mela, PetMel, JdlR} and the references they contain. In this paper we study a specific class of Bratteli-Vershik transformations known as \emph{limited scope adic transformations} which have zero entropy and contains all of the transformations contained in \cite{Mela}, including the Pascal adic transformation.  We show that the dimension group associated to the Bratteli diagram on which the limited scope adic are defined is order isomorphic to the continuous functions from the infinite path space into the integers modulo the continuous coboundaries.  This is a extension of the result of Herman, Putnam, and Skau for minimal Cantor systems, \cite{HPS}.  We also compute the dimension group directly for a particular subclass of limited scope adic transformations which are associated to polynomials over the natural numbers, defined below.  The transformations contained in \cite{Mela}, including the Pascal adic, are contained in this subclass.  Different limited scope adic transformations may have very different dynamical properties, and we establish several dynamical properties for the subclass of limited scope adic transformations associated to polynomials over the natural numbers as well as for the Euler adic transformation, which is defined below.  In particular we determine all of the ergodic adic invariant measures for the subclass of limited scope adic transformations associated to polynomials over the natural numbers.  We show that many of the adic transformations associated to polynomials over the natural numbers have roots of unity as eigenvalues where as the Euler adic transformation is totally ergodic.  We conclude by showing that all of the adic transformations associated to polynomials over the natural numbers as well as the Euler adic transformation is loosely Bernoulli.  This paper is based on the Ph.D. dissertation of the author at the University of North Carolina at Chapel Hill under the supervision of Karl Petersen, \cite{SBF1}.  

 Let $\VE$ be a Bratteli diagram such that for a
constant $d\in \N$, the number of
vertices at level $n$ is $nd+1$ and each vertex, labeled $(n,k)$,
$0\leq k\leq dn$, is connected by some positive number of edges to
each vertex in $(n+1,k+i)$ for all $i\in \{0,1,2,\dots,d\}$, and
there are no edges elsewhere.  We denote this family of Bratteli
diagrams by $\D$.

\begin{figure}[h]
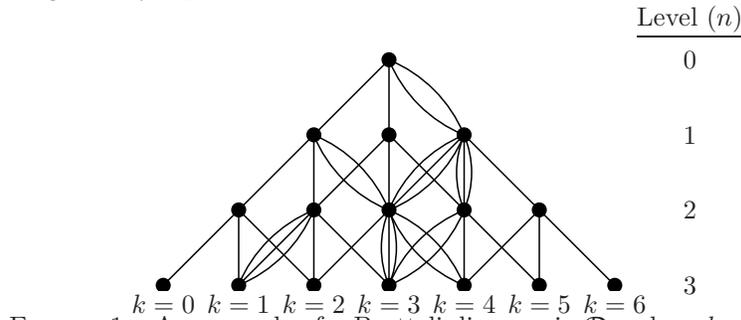

\begin{center}
\begin{graph}(6,3)(0,.25)
\unitlength=\unitlength
\roundnode{1}(3,3)\roundnode{2}(2,2)\roundnode{3}(3,2)\roundnode{4}(4,2)
\roundnode{5}(1,1)\roundnode{6}(2,1)\roundnode{7}(3,1)\roundnode{8}(4,1)\roundnode{9}(5,1)
\roundnode{10}(0,0)\roundnode{11}(1,0)\roundnode{12}(2,0)\roundnode{13}(3,0)
\roundnode{14}(4,0)\roundnode{15}(5,0)\roundnode{16}(6,0)
\edge{1}{2}\edge{1}{3}\bow{1}{4}{-.1}\bow{1}{4}{.1}\edge{2}{5}\edge{2}{6}
\bow{2}{7}{-.1}\bow{2}{7}{.1}\edge{3}{6}\edge{3}{7}\edge{3}{8}\edge{4}{7}
\bow{4}{7}{-.1}\bow{4}{7}{.1}\bow{4}{8}{-.1}\bow{4}{8}{.1}\edge{4}{8}\edge{4}{9}\edge{5}{10}\edge{5}{11}\edge{5}{12}
\edge{6}{11}\bow{6}{11}{.1}\bow{6}{11}{-.1}\edge{6}{12}\edge{6}{13}\edge{7}{12}\edge{7}{13}
\bow{7}{13}{.1}\bow{7}{13}{-.1}\bow{7}{14}{-.1}\bow{7}{14}{.1}\bow{8}{13}{-.1}\bow{8}{13}{.1}\edge{8}{14}\edge{8}{15}
\edge{9}{14}\edge{9}{15}\edge{9}{16}
\nodetext{10}(0,-.25){$k=0$}\nodetext{11}(0,-.25){$k=1$}\nodetext{12}(0,-.25){$k=2$}\nodetext{13}(0,-.25){$k=3$}
\nodetext{14}(0,-.25){$k=4$}\nodetext{15}(0,-.25){$k=5$}\nodetext{16}(0,-.25){$k=6$}
\freetext(7,3.5){\underline{Level
$(n)$}}\freetext(7,3){0}\freetext(7,2){1}\freetext(7,1){2}\freetext(7,0){3}
\end{graph}
\end{center}
\caption[A Bratteli diagram in $\D$]{ An example of a Bratteli
diagram in $\D$ when $d=2$.}
\end{figure}

For any diagram $\VE\in\D$, $X$ is the space of infinite edge paths
on $\VE$.   The vertex through which $\gamma$ passes at level $n$ is
denoted $(n,k_n(\gamma))$.  $X$ is a metric space with the standard
metric:  for $\gamma=\gamma_0\gamma_1\dots$ and $\xi=\xi_0\xi_1\dots$,
$d(\gamma, \xi)=2^{-j}$, where $j=\inf\{\gamma_j\neq
\xi_j\}$.

  An ordering given to edges of the diagram which
terminate into the same vertex is extended to a partial ordering on
the entire path space. Two paths $\gamma$ and $\xi$ are
\emph{comparable} if they agree after some level $n$ and disagree on
level $n-1$.  We then define $\gamma<\xi$ if and only if
$\gamma_{n-1}<\xi_{n-1}$ with respect to the edge ordering.  When we have endowed a Bratteli diagram $\VE$ with an edge ordering extending to infinite paths, we say that $\VE$ is an ordered Bratteli diagram and denote it by $\OVE$.   The
diagrams in $\D$ are drawn so that edges with the same range
increase in order from left to right.

 We denote by $X_{\max}$ the set of paths in $X$ for which
 all edges are maximal with respect to their edge ordering.
 Likewise we denote by $X_{\min}$ the set of paths in $X$ for which all edges
 are minimal with respect to their edge ordering.  For Bratteli diagrams in $\D$
 there are a countable number of paths in
  $X_{\max}\cup X_{\text{min}}$.  Indeed, for every $k$ in the
  set $\{0,1,\dots\}\cup\{\infty\}$ there is a unique associated
  path in $X_{\max}$, denoted $\gamma^k_{\max}$, which is defined as follows.
  For $k\neq \infty$ $\gamma_{\max}^k$ is the
  path in $X$ that travels down the far right side of the graph,
  following maximal edges, to level $n_0-1$, where $n_0\in \N$ is such
  that $(n_0-1)d<k\leq n_0d$, and then connects to vertex
  $(n_0,k)$ along the maximal edge.  Then for $n\geq n_0$,
  $k_n(\gamma_{\max}^k)=k$, and $\gamma_{\max}^k$ follows a maximal edge.  The path
  $\gamma_{\max}^{\infty}$ is the path which travels through the vertices
  $(n,dn)$ along maximal edges for all $n\in \N$.  Likewise for every $k$ in the set
  $\{0,1,\dots\}\cup\{\infty\}$ there is a unique path in
  $X_{\text{min}}$ denoted $\gamma^k_{\text{min}}$.  For $k\neq \infty$ this is the path in
  $X$ that travels down the left side of the graph along minimal
  edges to level $n_0-1$, where $n_0\in \N$ is such that
  $(n_0-1)d<k\leq n_0d$, and then connects to vertex
  $(n_0,n_0d-k)$ along the minimal edge.  Then for $n\geq n_0$,
  $k_n(\gamma^k_{\text{min}})=nd-k$.  The path $\gamma^\infty_{\min}$ is the path which
  travels through the vertices $(n,0)$ along minimal edges for all
  $n\in \N$.

\begin{figure}[h]
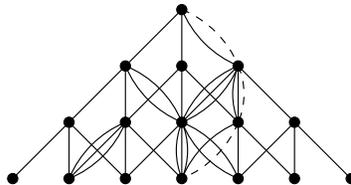

\begin{center}
\begin{graph}(4,2)(0,.5)
\unitlength=.75\unitlength
\roundnode{1}(3,3)\roundnode{2}(2,2)\roundnode{3}(3,2)\roundnode{4}(4,2)
\roundnode{5}(1,1)\roundnode{6}(2,1)\roundnode{7}(3,1)\roundnode{8}(4,1)\roundnode{9}(5,1)
\roundnode{10}(0,0)\roundnode{11}(1,0)\roundnode{12}(2,0)\roundnode{13}(3,0)
\roundnode{14}(4,0)\roundnode{15}(5,0)\roundnode{16}(6,0)
\edge{1}{2}\edge{1}{3}\bow{1}{4}{-.1}\bow{1}{4}{.1}[\graphlinedash{
3 3}]\edge{2}{5}\edge{2}{6}
\bow{2}{7}{-.1}\bow{2}{7}{.1}\edge{3}{6}\edge{3}{7}\edge{3}{8}\edge{4}{7}
\bow{4}{7}{-.1}\bow{4}{7}{.1}\bow{4}{8}{-.1}\bow{4}{8}{.1}[\graphlinedash{
3 3}]\edge{4}{8}\edge{4}{9}\edge{5}{10}\edge{5}{11}\edge{5}{12}
\edge{6}{11}\bow{6}{11}{.1}\bow{6}{11}{-.1}\edge{6}{12}\edge{6}{13}\edge{7}{12}\edge{7}{13}
\bow{7}{13}{.1}\bow{7}{13}{-.1}\bow{7}{14}{-.1}\bow{7}{14}{.1}\bow{8}{13}{-.1}\bow{8}{13}{.1}[\graphlinedash{
3 3}]\edge{8}{14}\edge{8}{15} \edge{9}{14}\edge{9}{15}\edge{9}{16}
\end{graph}
\end{center}
\caption{ The
dashed path is the first three edges of $\gamma_{\max}^3$.}
\end{figure}

Let $T:X\to X$ be the Bratteli-Vershik transformation
which maps a path $\gamma$ to the next largest path with respect to
the partial ordering of edges, if one exists.  If such a path exits,
we will call it the \emph{successor} of $\gamma$.  All paths in
$X\setminus X_{\max}$ have unique successors, and hence $T$ is well
defined off of $X_{\max}$.  We will now define $T$ on $X_{\max}$ so
that
  $T(\gamma^k_{\max})=\gamma^k_{\min}$ for $0<k<\infty$,
  $T(\gamma^0_{\max})=\gamma^\infty_{\text{min}}$, and
  $T(\gamma^\infty_{\max})=\gamma^0_{\text{min}}$.  In this way $T$ is a
  bijection on the whole space $X$; but not continuous on
  $X_{\max}$.  The family of Bratteli-Vershik systems determined in the manner from
  Bratteli diagrams in $\D$ are said to be of limited scope and will be denoted $\Sl$.

For $(X,T)\in \Sl$, we say a path $\gamma\in X$ is \emph{eventually
diagonal to the left} if there exists and $N\geq 0$ such that for
$n\geq N$, $k_n(\gamma)=k_N(\gamma)$.  We say a path $y\in X$ is
\emph{eventually diagonal to the right} if there exists an $M\geq 0$
such that for $m\geq M$,
$k_m(\gamma)=dm-k_M(\gamma)$.\label{evendiagdef}  We will say that a
path is \emph{eventually diagonal} if the direction is either clear
or unknown.  All paths in the orbits of $X_{\max}$ and $X_{\min}$
are eventually diagonal.

\begin{prop} For every $\gamma\in X$, exactly one of
the following holds.\\
1. $\gamma$ is eventually diagonal to the right.\\
2. $\gamma$ is eventually diagonal to the left.\\
3. $\overline{\mathcal{O}(\gamma)}=X$\label{prop311}
\end{prop}
\begin{proof}
  If $\gamma$ is not eventually diagonal, both $k_n(\gamma)$ and $dn-k_n(\gamma)$ are unbounded.
  Then for any $\xi\in X$ and $m\in\N$ there is an $n_0>m$ such that
  $k_m(\xi)\leq k_{n_0}(\gamma)$ and $dm-k_m(\xi)\leq dn_0-k_{n_0}(\gamma)$.  Hence,
  $k_{n_0}(\gamma)-d(n_0-m)\leq k_m(\xi)\leq k_{n_0}(\gamma)$.  Therefore
  there is a path from $(m,k_m(\xi))$ to $(n_0,k_{n_0}(\gamma))$.  Then there is
  a $j\in \Z$ so that $T^j\gamma$ coincides with $\xi$ along the first
  $m$ edges, showing that $\mathcal{O}(\gamma)$ is dense in $X$.

If $\gamma$ is eventually diagonal to the right (resp. to the
left), there exists an $N\in\N$ such that for any $\eta\in
\mathcal{O}(\gamma)$ and all $n\in \N$, $k_n(\eta)<N$ or
$k_n(\eta)>dn-N$. Now choose $\xi\in X$ and $m\in \N$ for which
$N<k_m(\xi)<dm-N$ and let $B_{2^{-m}}(\xi)$ be the ball of radius
$2^{-m}$ around $\xi$.  Then $\mathcal{O}(\gamma)\cap
B_{2^{-m}}(\xi)=\emptyset$. Hence, $\mathcal{O}(\gamma)$ is not
dense in $X$.
\end{proof}

\section{Examples}
  In this section we give a description of some examples of
  Bratteli-diagrams in $\D$.  Endowed with the above adic transformation
  they generate Bratteli-Vershik systems in $\Sl$.  Specifically we give
  examples of those Bratteli-Vershik systems determined
  by polynomials over $\N$, the Euler adic, and the reverse Euler
  adic.

  We begin with those determined by polynomials over $\N$.
Every positive integer polynomial of degree $d$ determines a
Bratteli-Vershik system in $\Sl$; we will denote this subfamily of
systems by $\Slp$ and the diagrams by $\Dlp$.

Let $a_0,a_1,...,a_d\in \N$ and $p(x)=a_0+a_1x+...+a_dx^d$. The
\emph{Bratteli diagram associated to $p(x)$} is a Bratteli diagram
in $\D$ such that for every level $n$, the number of vertices is
$dn+1$ and the number of edges from $(n,k)$ to $(n+1,k+j)$ is $a_j$,
with $a_j$=0 for $j>d$ and $j<0$.

\begin{figure}[h]
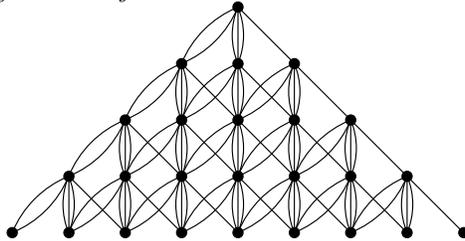

\begin{center}
\begin{graph}(9,2)(-1.5,.5)
\unitlength=.75\unitlength
\roundnode{0}(4,4)\roundnode{1}(3,3)\roundnode{2}(4,3)\roundnode{3}(5,3)
\roundnode{4}(2,2)\roundnode{5}(3,2)\roundnode{6}(4,2)\roundnode{7}(5,2)
\roundnode{8}(6,2)\roundnode{9}(1,1)\roundnode{10}(2,1)\roundnode{11}(3,1)\roundnode{12}(4,1)
\roundnode{13}(5,1)\roundnode{14}(6,1)\roundnode{15}(7,1)
\roundnode{16}(0,0)\roundnode{17}(1,0)\roundnode{18}(2,0)\roundnode{19}(3,0)
\roundnode{20}(4,0)\roundnode{21}(5,0)\roundnode{22}(6,0)\roundnode{23}(7,0)\roundnode{24}(8,0)
\bow{0}{1}{-.1}\bow{0}{1}{.1}\bow{0}{2}{.1}\edge{0}{2}\bow{0}{2}{-.1}\edge{0}{3}
\bow{1}{4}{-.1}\bow{1}{4}{.1}\bow{1}{5}{.1}\edge{1}{5}\bow{1}{5}{-.1}\edge{1}{6}
\bow{2}{5}{-.1}\bow{2}{5}{.1}\bow{2}{6}{.1}\edge{2}{6}\bow{2}{6}{-.1}\edge{2}{7}
\bow{3}{6}{-.1}\bow{3}{6}{.1}\bow{3}{7}{.1}\edge{3}{7}\bow{3}{7}{-.1}\edge{3}{8}
\bow{4}{9}{-.1}\bow{4}{9}{.1}\bow{4}{10}{.1}\edge{4}{10}\bow{4}{10}{-.1}\edge{4}{11}
\bow{5}{10}{-.1}\bow{5}{10}{.1}\bow{5}{11}{.1}\edge{5}{11}\bow{5}{11}{-.1}\edge{5}{12}
\bow{6}{11}{-.1}\bow{6}{11}{.1}\bow{6}{12}{.1}\edge{6}{12}\bow{6}{12}{-.1}\edge{6}{13}
\bow{7}{12}{-.1}\bow{7}{12}{.1}\bow{7}{13}{.1}\edge{7}{13}\bow{7}{13}{-.1}\edge{7}{14}
\bow{8}{13}{-.1}\bow{8}{13}{.1}\bow{8}{14}{.1}\edge{8}{14}\bow{8}{14}{-.1}\edge{8}{15}
\bow{9}{16}{-.1}\bow{9}{16}{.1}\bow{9}{17}{.1}\edge{9}{17}\bow{9}{17}{-.1}\edge{9}{18}
\bow{10}{17}{-.1}\bow{10}{17}{.1}\bow{10}{18}{.1}\edge{10}{18}\bow{10}{18}{-.1}\edge{10}{19}
\bow{11}{18}{-.1}\bow{11}{18}{.1}\bow{11}{19}{.1}\edge{11}{19}\bow{11}{19}{-.1}\edge{11}{20}
\bow{12}{19}{-.1}\bow{12}{19}{.1}\bow{12}{20}{.1}\edge{12}{20}\bow{12}{20}{-.1}\edge{12}{21}
\bow{13}{20}{-.1}\bow{13}{20}{.1}\bow{13}{21}{.1}\edge{13}{21}\bow{13}{21}{-.1}\edge{13}{22}
\bow{14}{21}{-.1}\bow{14}{21}{.1}\bow{14}{22}{.1}\edge{14}{22}\bow{14}{22}{-.1}\edge{14}{23}
\bow{15}{22}{-.1}\bow{15}{22}{.1}\bow{15}{23}{.1}\edge{15}{23}\bow{15}{23}{-.1}\edge{15}{24}
\end{graph}
\end{center}
\caption{The first five levels of $\VE_{2+3x+x^2}$}
\end{figure}

These diagrams have the property that for any vertex $(n,k)$ the
number of paths from the root vertex, (0,0), into $(n,k)$ is the coefficient
of $x^k$ in the polynomial $\left(p(x)\right)^n$.  The most famous example of this
is the Pascal adic.

Our next example is the Euler adic, introduced in \cite{BKPS}.  The
Euler graph, is the Bratteli diagram in $\D$ for which the number of
vertices at each level $n$ is $n+1$, and the number of edges
connecting vertex $(n,k)$ to vertex $(n+1,k)$ is $k+1$ while the
number of edges connecting vertex $(n,k)$ to vertex $(n+1,k+1)$ is
$n-k+1$.

\begin{figure}[h]
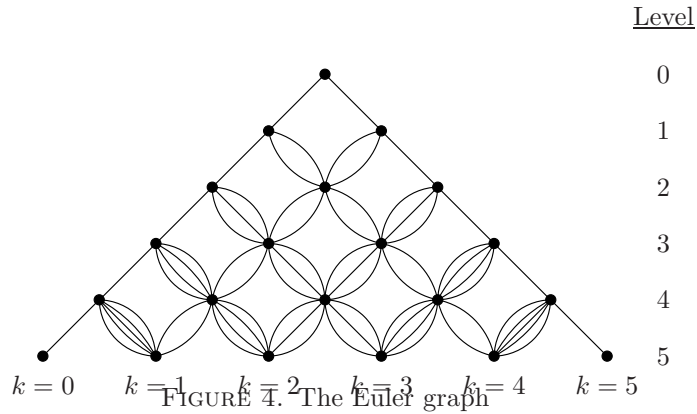

\begin{center}
\begin{graph}(6,4)(0,-.5)
\unitlength=.75\unitlength
 \roundnode{v00}(4,4)
\roundnode{v10}(3,3)\roundnode{v11}(5,3)
\roundnode{v20}(2,2)\roundnode{v21}(4,2)\roundnode{v22}(6,2)\edge{v00}{v10}\edge{v00}{v11}
\edge{v10}{v20}\bow{v10}{v21}{.15}
\bow{v10}{v21}{-.15}\bow{v11}{v21}{-.15}
\bow{v11}{v21}{.15}\edge{v11}{v22}\roundnode{v30}(1,1)\roundnode{v31}(3,1)\roundnode{v32}(5,1)\roundnode{v33}(7,1)
\edge{v20}{v30}\edge{v20}{v31}\bow{v20}{v31}{.2}\bow{v20}{v31}{-.2}\bow{v21}{v31}{.15}\bow{v21}{v31}{-.15}
\bow{v21}{v32}{.15}\bow{v21}{v32}{-.15}\edge{v22}{v32}\bow{v22}{v32}{.2}\bow{v22}{v32}{-.2}\edge{v22}{v33}
\roundnode{v40}(0,0)\roundnode{v41}(2,0)\roundnode{v42}(4,0)\roundnode{v43}(6,0)\roundnode{v44}(8,0)
\edge{v30}{v40}\bow{v30}{v41}{.075}\bow{v30}{v41}{.2}\bow{v30}{v41}{-.075}\bow{v30}{v41}{-.2}\bow{v31}{v41}{.15}
\bow{v31}{v41}{-.15}\bow{v31}{v42}{.2}\bow{v31}{v42}{-.2}\edge{v31}{v42}\edge{v32}{v42}\bow{v32}{v42}{.2}
\bow{v32}{v42}{-.2} \bow{v32}{v43}{.15}\bow{v32}{v43}{-.15}
\bow{v33}{v43}{.075}\bow{v33}{v43}{.2}\bow{v33}{v43}{-.075}\bow{v33}{v43}{-.2}\edge{v33}{v44}
\roundnode{v50}(-1,-1)\roundnode{v51}(1,-1)\roundnode{v52}(3,-1)\roundnode{v53}(5,-1)\roundnode{v54}(7,-1)
\roundnode{v55}(9,-1)
\edge{v40}{v50}\edge{v40}{v51}\bow{v40}{v51}{.075}\bow{v40}{v51}{.2}\bow{v40}{v51}{-.075}\bow{v40}{v51}{-.2}
\bow{v41}{v51}{.15}\bow{v41}{v51}{-.15}
\bow{v41}{v52}{.075}\bow{v41}{v52}{.2}\bow{v41}{v52}{-.075}\bow{v41}{v52}{-.2}
\edge{v42}{v52}\bow{v42}{v52}{.2}\bow{v42}{v52}{-.2}
\edge{v42}{v53}\bow{v42}{v53}{.2}\bow{v42}{v53}{-.2}
\bow{v43}{v53}{.075}\bow{v43}{v53}{.2}\bow{v43}{v53}{-.075}\bow{v43}{v53}{-.2}
\bow{v43}{v54}{.15}\bow{v43}{v54}{-.15}
\edge{v44}{v54}\bow{v44}{v54}{.075}\bow{v44}{v54}{.2}\bow{v44}{v54}{-.075}\bow{v44}{v54}{-.2}
\edge{v44}{v55} \freetext(10,5){\underline{Level}}\freetext(10,4){0}
\freetext(10,3){1}\freetext(10,2){2}\freetext(10,1){3}\freetext(10,0){4}\freetext(10,-1){5}
\nodetext{v50}(0,-.5){$k=0$}
\nodetext{v51}(0,-.5){$k=1$}\nodetext{v52}(0,-.50){$k=2$}\nodetext{v53}(0,-.50){$k=3$}
\nodetext{v54}(0,-.50){$k=4$}\nodetext{v55}(0,-.50){$k=5$}
\end{graph}
\end{center}
\caption{The Euler graph} \label{Graph}
\end{figure}

The Euler graph has the property that the number of paths from the
root vertex to a vertex $(n,k)$ is the \emph{Eulerian number}
$A(n,k)$.  That is, the number of permutations $i_1i_2\dots i_{n+1}$
of $\{1, 2, \dots, n+1\}$ with exactly $k$ rises and $n-k$ falls;
see \cite{Comtet} for background concerning Eulerian numbers.

The final example is the reverse Euler adic. The reverse Euler
graph, is the Bratteli diagram in $\D$ for which the number of
vertices at each level $n$ is $n+1$, and the number of edges
connecting vertex $(n,k)$ to vertex $(n+1,k)$ is $n-k+1$ while the
number of edges connecting vertex $(n,k)$ to vertex $(n+1,k+1)$ is
$k+1$; the reverse of the connection in the Euler graph.

\begin{figure}[h]
\begin{center}
\begin{graph}(6,4)(0,-.5)
\unitlength=.75\unitlength
 \roundnode{v00}(4,4)
\roundnode{v10}(3,3)\roundnode{v11}(5,3)
\roundnode{v20}(2,2)\roundnode{v21}(4,2)\roundnode{v22}(6,2)\edge{v00}{v10}\edge{v00}{v11}
\edge{v10}{v21}\bow{v10}{v20}{.15}
\bow{v10}{v20}{-.15}\bow{v11}{v22}{-.15}
\bow{v11}{v22}{.15}\edge{v11}{v21}\roundnode{v30}(1,1)\roundnode{v31}(3,1)\roundnode{v32}(5,1)\roundnode{v33}(7,1)
\edge{v20}{v31}\edge{v20}{v30}\bow{v20}{v30}{.2}\bow{v20}{v30}{-.2}\bow{v21}{v32}{.15}\bow{v21}{v32}{-.15}
\bow{v21}{v31}{.15}\bow{v21}{v31}{-.15}\edge{v22}{v33}\bow{v22}{v33}{.2}\bow{v22}{v33}{-.2}\edge{v22}{v32}
\roundnode{v40}(0,0)\roundnode{v41}(2,0)\roundnode{v42}(4,0)\roundnode{v43}(6,0)\roundnode{v44}(8,0)
\edge{v30}{v41}\bow{v30}{v40}{.075}\bow{v30}{v40}{.2}\bow{v30}{v40}{-.075}\bow{v30}{v40}{-.2}\bow{v31}{v42}{.15}
\bow{v31}{v42}{-.15}\bow{v31}{v41}{.2}\bow{v31}{v41}{-.2}\edge{v31}{v41}\edge{v32}{v43}\bow{v32}{v43}{.2}
\bow{v32}{v43}{-.2} \bow{v32}{v42}{.15}\bow{v32}{v42}{-.15}
\bow{v33}{v44}{.075}\bow{v33}{v44}{.2}\bow{v33}{v44}{-.075}\bow{v33}{v44}{-.2}\edge{v33}{v43}
\roundnode{v50}(-1,-1)\roundnode{v51}(1,-1)\roundnode{v52}(3,-1)\roundnode{v53}(5,-1)\roundnode{v54}(7,-1)
\roundnode{v55}(9,-1)
\edge{v40}{v51}\edge{v40}{v50}\bow{v40}{v50}{.075}\bow{v40}{v50}{.2}\bow{v40}{v50}{-.075}\bow{v40}{v50}{-.2}
\bow{v41}{v52}{.15}\bow{v41}{v52}{-.15}
\bow{v41}{v51}{.075}\bow{v41}{v51}{.2}\bow{v41}{v51}{-.075}\bow{v41}{v51}{-.2}
\edge{v42}{v53}\bow{v42}{v53}{.2}\bow{v42}{v53}{-.2}
\edge{v42}{v52}\bow{v42}{v52}{.2}\bow{v42}{v52}{-.2}
\bow{v43}{v54}{.075}\bow{v43}{v54}{.2}\bow{v43}{v54}{-.075}\bow{v43}{v54}{-.2}
\bow{v43}{v53}{.15}\bow{v43}{v53}{-.15}
\edge{v44}{v55}\bow{v44}{v55}{.075}\bow{v44}{v55}{.2}\bow{v44}{v55}{-.075}\bow{v44}{v55}{-.2}
\edge{v44}{v54} \freetext(10,5){\underline{Level}}\freetext(10,4){0}
\freetext(10,3){1}\freetext(10,2){2}\freetext(10,1){3}\freetext(10,0){4}\freetext(10,-1){5}
\nodetext{v50}(0,-.5){$k=0$}
\nodetext{v51}(0,-.5){$k=1$}\nodetext{v52}(0,-.50){$k=2$}\nodetext{v53}(0,-.50){$k=3$}
\nodetext{v54}(0,-.50){$k=4$}\nodetext{v55}(0,-.50){$k=5$}
\end{graph}
\end{center}
\caption{The Reverse Euler graph} \label{reversed}
\end{figure}

The reverse Euler graph has the property that the number of paths
from the root vertex to a vertex $(n,k)$ is $(n+1)!$.

\section{Dimension Groups}
\label{dgsection}
Every Bratteli diagram can be completely described by a sequence of
\emph{incidence} matrices.  For any pair of consecutive levels $n-1$
and $n$, with $v_{n-1}$ vertices on level $n-1$ and $v_{n}$ on level
$n$, the incidence matrix $D_n$ is a $v_n\times v_{n-1}$ matrix such
that $[D_n]_{i,j}$ is the number of edges connecting vertices
$(n-1,j)$ and $(n,i)$.

For every Bratteli diagram $\VE$ there is an associated ordered group called the dimension
group and denoted by $K_0\VE$.  Explicitly it is the direct limit of the following
directed system:
$$\ds \Z^{|V_0|=1}\xrightarrow{\phi_1}
\Z^{|V_1|}\xrightarrow{\phi_2}\Z^{|V_2|}\xrightarrow{\phi_3}...$$
where  for each $i=1,2,\dots $ $\phi_i$ is the group homomorphism
determined by the incidence matrix between levels $i-1$ and $i$ of
the Bratteli diagram.  The positive set consists of the equivalence
classes for which there is a nonnegative vector representative.  The
equivalence class of $1\in \Z$ is called the \emph{distinguished
order unit}.  It is of interest to note that $K_0\VE$ is not
dependent on the ordering of the edges or on the associated
dynamical system.  For further references on
ordered groups and dimension groups see  \cite{HPS},\cite{GPS}, \cite{DHS} and \cite{Boyle}.

 If $(X,T)$ is any dynamical system
let $C(X,\Z)$ denote the additive group of continuous functions from
the space $X$ to $\Z$ and define
$$\partial_T C(X,\Z)=\{g\circ T-g| g\in
  C(X,\Z)\}.$$  The elements of $\partial_T C(X,\Z)$ are called the \emph{coboundaries} of $(X,\phi)$.
In the case that $T$
  is a homeomorphism,  $K^0(X,T)$ is defined
  to be $C(X,\Z)/\partial_T C(X,\Z)$.

If an ordered Bratteli diagram has $X_{\max}$ and $X_{\min}$ both
one point sets, then the diagram is said to be \emph{essentially
simple}.

  \begin{theo}[Herman, Putnam, Skau \cite{HPS}]  Let $(\V,\E,\geq)$ be an essentially simple
  ordered Bratteli diagram and let $(X,T)$ be its associated
  Bratteli-Vershik system.  Then there is an order isomorphism
  $$\theta:K_0\VE\to K^0(X,\phi)$$
  which maps the distinguished order unit of
  $K_0\VE$ to the equivalence class of the constant function 1.
\label{dgt}  \end{theo}

More discussion of this result can be found in \cite{GPS, GW}.

 In the case of systems in $\Sl$, which are not essentially
simple, $\partial_{T}C(X,\Z)$ may not be not contained in $C(X,\Z)$
as $T$ is not continuous everywhere. Nevertheless, by slightly
adjusting the definition of $K^0(X,T)$ to be
$C(X,\Z)/(\partial_{T}C(X,\Z)\cap C(X,\Z))$ we can achieve a result
similar to Theorem \ref{dgt}.

\begin{theo}
 For $(X,T) \in\Sl$, there is an order isomorphism
$$K_0(\mathcal{V,E}) \cong K^0(X,T)$$
which maps the distinguished order unit of $K_0\VE$ to the
equivalence class of the constant function 1.\label{dgtheorem1}
\end{theo}

Before we give the proof, we introduce some useful notation.  Let $\OVE$ be an ordered Bratteli diagram.
For $n=1,2,\dots$ and $0\leq k\leq dn$ define $\dim(n,k)$ to be the
number of finite paths from the root vertex, (0,0) into the vertex $(n,k)$.  For any
vertex $(n,k)\in \V$ there is a cylinder determined by the path from
the root vertex to $(n,k)$ for which all the edges are minimal (maximal).  We will call this the \emph{minimal
(maximal) cylinder terminating at vertex $(n,k)$}. Denote by
$Y_n(k,0)$ the minimal cylinder into vertex $(n,k)$, and let
$Y_n(k,i)=T^i(Y_n(k,0))$ for $i=0,1,\dots, \dim(n,k)-1$. For each
$n=0,1,2,\dots,$ denote the union of all the minimal cylinders of
length $n$ by $Y_n$, so that
$$\ds Y_n=\bigcup_{0\leq k\leq |\V_n|-1}Y_n(k,0).$$

\begin{proof} This proof is an adaptation of the dynamical proof of Theorem
\ref{dgt} given by Glasner and Weiss in \cite{GW}. We will first
define a group homomorphism $J:C(X,\Z)\to K_0\VE$. Then we will
define a set $B$ and show that it is a subset of $C(X,\Z)$. Then we
will show $B=\ker(J)$ by first showing $B\subset \ker(J)$ and then
$\ker(J)\subset B$. This will induce a one-to-one group homomorphism
$\tilde J:C(X,\Z)/B\to K_0\VE$.  We will then show that $\tilde J$
is surjective and in fact an order isomorphism. Lastly we will show
$B=\partial_TC(X,\Z)\cap C(X,\Z)$.

  Let $f\in C(X,\Z)$.  Since $X$ is compact, $f$ is bounded and
hence takes on only finitely many values.  Let $\{l_1,\dots l_j\}$
be the set of these values and let $U_i=f^{-1}\{l_i\}$ for each $i$.
If $i=1,\dots,j$ and $\gamma\in U_i$, then there is a cylinder set
$C_\gamma\subset U_i$ of the form $[c_0c_1\dots c_{N_\gamma-1}]$
which contains $\gamma$.  From $\{C_\gamma|\gamma\in X\}$ select a
finite subcover $\{C_{\gamma^1},C_{\gamma^2},\dots,C_{\gamma^r}\}$.
Then for some $i\in\{1,2,\dots,r\}$, $C_{\gamma^i}$ is of longest
length, $N_1(f)$, and $f$ is constant on any cylinder of length
$n\geq N_1(f)$.  For $n\geq N_1(f)$ define an element $\tilde f_n
\in \Z^{dn+1}$ by letting, for each $0\leq k\leq dn$ and any
$\gamma\in Y_n(k,0)$,
$$\tilde f_n(k)=f(\gamma) + f(T\gamma)+f(T^2\gamma)+\dots + f(T^{\dim(n,k)-1}\gamma).$$

$D_n$ denotes the adjacency matrix of the edges
connecting levels $n-1$ and $n$.  Then
$$\tilde f_{n+1}(i)=\sum_{j=0}^{nd}\tilde
f_n(j)(D_n)_{i,j}=(\tilde f_n D_n)(i).$$ Therefore the sequence
$\tilde f_n$ defines an element $J(f)\in K_0\VE$. Clearly $J:C(X,\Z)\to K_0\VE$ is a group homomorphism.

Let $G=\{g\in C(X,\Z)|\exists\
N_2(g)\text{ such that }n\geq N_2(g) \implies\forall \gamma\in
Y_n,\ g(\gamma)=c\}.$  In other words, $g$ takes the same value on
all the minimal cylinders into level $n$.  Define $B=\{g\circ
T-g|g\in
 G\}$.

 We now show that $B\subset C(X,\Z)$.  For $f\in B$ with $f=g\circ
T-g$, $f$ is continuous on $X\setminus X_{\max}$, we
need to check continuity of $f$ on $X_{\max}$.  Let $m\geq
\max\{N_1(g),N_2(g)\}$ be such that $g$ is constant on each cylinder
of length $m$ and $g$ is also constant on $Y_m$. For
$\gamma_{\max}\in X_{\max}$ and $\xi\in X$,
$d(\gamma_{\max},\xi)<2^{-m}$ implies that $\gamma_{\max}$ and $\xi$
are both in the same maximal cylinder terminating at vertex
$(m,k_m(\gamma_{\max}))$, and hence $g(\gamma_{\max})=g(\xi)$. Since
$T(\gamma_{\max})$ and $T(\xi)$ are both in $Y_m$, we have $(g\circ
T)(\gamma)=(g\circ T)(\xi)$. Hence $f(\gamma)=f(\xi)$, and so $f$ is
continuous.

 We will show that $B=\ker(J)$.  If $f=g\circ T-g\in B$,
$n\geq \max\{N_1(g)N_2(g)\}$, $0\leq k\leq dn$, and any $\gamma\in
Y_n(k,0)$, then
$$\tilde f_n(k)=f(\gamma) + f(T\gamma)+f(T^2\gamma)+\dots + f(T^{\dim(n,k)-1}\gamma)=g\circ
T^{\dim(n,k)}(\gamma)-g(\gamma).$$ Since both $\gamma$ and
$T^{\dim(n,k)}(\gamma)\in Y_n$ and $g$ is constant on $Y_n$, $\tilde
f_n(k)=0$. Therefore $J(f)=0$, which implies $B\subset \ker(J)$.

Conversely, if $f\in C(X,\Z)$ and $J(f)=0$, there is an $n>N_1(f)$
for which $\tilde f_n=0$. We will define a function $g\in C(X,\Z)$
so that $f=g\circ T-g$. Let $g=0$ on $Y_n$.  For $1\leq l\leq
\dim(n,k)$, choose any $\gamma\in Y_n(k,0)$ and let $g\equiv
f(\gamma)+f(T\gamma)+\dots+f(T^{l-1}\gamma)$ on $Y_n(k,l)$. Now $g$
is everywhere defined, and clearly $f=g\circ T-g$ on every cylinder
terminating at vertex $(n,k)$ except maybe on the maximal cylinder.
However, for $\gamma\in Y_n(k,0)$, $g(T^{\dim(n,k)}\gamma)=0$ and
$\tilde f_n(k)=0$, so we have
\[\begin{array}{lll}
&&g(T^{\dim(n,k)}\gamma)-g(T^{\dim(n,k)-1}\gamma)\\&=&-g(T^{\dim(n,k)-1}\gamma)\\
&=&-(f(\gamma)+f(T\gamma)+\dots +f(T^{\dim(n,k)-2}\gamma))\\
&=&-(f(\gamma)+f(T\gamma)+\dots +f(T^{\dim(n,k)-1}\gamma))+f(T^{\dim(n,k)-1}\gamma)\\
&=&-\tilde f_n(k)+f(T^{\dim(n,k)-1}\gamma)\\
&=&f(T^{\dim(n,k)-1}\gamma).
\end{array}\]
Thus $f=g\circ T-g$ also on the maximal cylinder, and hence $f\in
B$. Thus $B= \ker(J)$, and $J$ induces an injective group
homomorphism $\tilde J: C(X,\Z)/B\to K_0\VE.$

We now show that $\tilde J$ is onto and an order isomorphism. Given
$a\in K_0\VE$, choose an $n\in \Z_+$ so that the equivalence class
$a$ has a representative $a_n\in\Z^{dn+1}$. Define $f$ as follows.
For $k=0,1,\dots,dn$ and $\gamma\in Y_n(k,0)$, let
$f(\gamma)=a_n(k)$ and elsewhere put $f=0$. Then $\tilde
f_n(k)=a_n(k)$, so that $J(f)=a$, and thus $\tilde J$ is onto.
Clearly $\tilde J$ takes positive elements to positive elements, and
the preceding argument shows that the unique preimage of every
positive element under $\tilde J$ is a positive element. Thus
$C(X,\Z)/B$ is order isomorphic to $K_0\VE$ by the map $\tilde J$,
which maps the equivalence class of the constant function 1 to the
distinguished order unit of $K_0\VE$.

\begin{figure}[h]
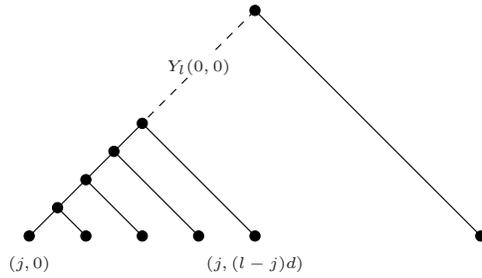

\begin{center}
\begin{graph}(6,2)(0,.25)
\unitlength=.75\unitlength
\roundnode{1}(4,4)\roundnode{2}(2,2)\roundnode{3}(0,0)\roundnode{4}(4,0)
\roundnode{5}(8,0)\roundnode{6}(1,0)\roundnode{7}(2,0)\roundnode{8}(3,0)
\roundnode{9}(.5,.5)\roundnode{10}(1,1)\roundnode{11}(1.5,1.5)\edge{9}{6}
\edge{10}{7}\edge{11}{8} \edge{1}{5}\edge{1}{2}[\graphlinedash{3 3
3}]\edgetext{1}{2}{\tiny$Y_l(0,0)$}
\edge{2}{3}\edge{2}{4}\nodetext{3}(0,-.5){\tiny$(j,0)$}\nodetext{4}(0,-.5){\tiny$(j,(l-j)d)$}
\end{graph}
\end{center}
\caption{ Connections from
level $l$ to $j$.} \label{minimal}
\end{figure}

We now show that $\partial_TC(X,\Z)\cap C(X,\Z)\subset B$. Let $f\in
\partial_T C(X,\Z)\cap C(X,\Z)$ be given. Then $f=g\circ T-g$ for some $g\in C(X,\Z)$, and $f$ is
continuous. We have to show that there is an $N_2(g)$ so that for
each $n\geq N_2(g)$, $g$ takes the same value on all of $Y_n$. Since
$g\in C(X,\Z)$, we can choose $l=N_1(g)$ such that $g$ is constant
on cylinder sets of length $l$. Then for every level $j\geq l$, and
every $i\in \{0,1,\dots, (j-l)d\}$, $Y_j(i,0)\subset Y_l(0,0)$, (see
Figure \ref{minimal}). Now consider $k<dl$ and  $\gamma_{\max}^k\in
X_{\max}$. Then $T(\gamma_{\max}^k)=\gamma_{\min}^k\in Y_l(ld-k,0)$.

Since $f=g\circ T-g\in C(X,\Z)$, given $\gamma_{\max}^k\in X$ with
$k<dl$ there is a $\delta>0$ such that
$d(\gamma_{\max}^k,\xi)<\delta$ implies $f(\gamma_{\max}^k)=f(\xi)$.
We will choose a $\xi\in X$ sufficiently close to $\gamma_{\max}^k$
such that $f(\gamma_{\max}^k)=f(\xi)$ and
$g(\gamma_{\max}^k)=g(\xi)$ which implies $g\circ
T(\gamma_{\max}^k=g\circ T(\xi)$.  Choose $j$ so that
$2^{-j}<\delta$ and $(j-l)d>k+1$. Now let $\xi$ be a path in $X$
such that $\xi_i=(\gamma_{\max}^k)_i$ for each $i=0,1,\dots,j-1$ and
$\xi_j\neq (\gamma_{\max}^k)_j$. Then
$d(\gamma_{\max}^k,\xi)<\delta$, so that
$f(\gamma_{\max}^k)=f(\xi)$. Since $j>l$, $\gamma_{\max}^k$ and
$\xi$ are in the same maximal cylinder which terminates at $(l,k)$,
which implies $g(\gamma_{\max}^k)=g(\xi)$.  Thus
$f(\gamma_{\max}^k)=f(\xi)$ implies $(g\circ
T)(\gamma_{\max}^k)=(g\circ T)(\xi)$. Since $s(\xi_j)=(j,k)$, and
$\xi_j$ is the first non-maximal edge of $\xi$, $T\xi$ is in either
$Y_j(k,0)$ or $Y_j(k+1,0)$ depending on the source of the successor
of $\xi_j$. Since $k+1<(j-l)d$, we have that
 $T\xi\in Y_l(0,0)$. Then $(g\circ
T)(\gamma_{\max}^k)=(g\circ T)(\xi)$, and $g$ constant on each
cylinder of length $l$ implies $g(Y_l(ld-k),0)=g(Y_l(0,0))$. Since
$k<dl$ was arbitrary, we have shown that $g$ is constant on all
$Y_l(k,0)$ for $k<dl$.  It remains only to show that $g$ takes this
same value on $Y_l(dl,0)$. Consider $\gamma_{\max}^\infty$, and
choose $j\geq l$ so that $2^{-j}<\delta$. Then
$d(\gamma_{\max}^{\infty},\gamma_{\max}^{jd})<\delta$, which implies
$\gamma_{\max}^{\infty}$ and $\gamma_{\max}^{jd}$ are both in the
maximal cylinder terminating at vertex $(l,dl)$. Thus
$g(\gamma_{\max}^{\infty})=g(\gamma_{\max}^{jd})$.  Then
$f(\gamma_{\max}^{\infty})=f(\gamma_{\max}^{jd})$ and
$g(\gamma_{\max}^{\infty})=g(\gamma_{\max}^{jd})$ implies $(g\circ
T)(\gamma_{\max}^{\infty})=(g\circ T)(\gamma_{\max}^{jd})$.  Thus
$T\gamma_{\max}^{\infty}\in Y_l(dl,0)$, $T\gamma_{\max}^{jd}\in
Y_l(0,0)$ and $g$ constant on each cylinder of length $l$ implies
$g(Y_l(0,0))=g(Y_l(dl,0))$.  Hence $g$ is constant on $Y_l$, as
required.
\end{proof}

We say that a Bratteli diagram is \emph{stationary} if for $n,m\geq 2$, the incidence matrices $D_n$ and $D_m$ are equal.  In this case there is a natural method of computing the dimension group.  See \cite{Boyle} for the exact construction.  In the case of  $\Dlp$ the Bratteli diagrams are clearly not stationary but in a certain sense the group homomorphisms that define the associated dimension groups are stationary and we have the following theorem.

\begin{theo}  The dimension group $K_0\vep$ associated to $\vep$ is order isomorphic to
the ordered group $\G$ of rational functions of the form
$$\ds \frac{r(x)}{p(x)^m},$$  where $r(x)$ is any polynomial with
integer coefficients such that $\deg(r(x))\leq md$. Addition of two
elements is given by $$\ds
\frac{r(x)}{p(x)^m}+\frac{s(x)}{p(x)^l}=\frac{r(x)+s(x)p(x)^{m-l}}{p(x)^m}$$
if $l\leq m$. The positive set $(\G)_+$ consists of the elements of
$\G$ such that there is an $l$ for which the numerator of
$$\ds\frac{r(x)(p(x))^l}{p(x)^{l+m}}$$ has all positive
coefficients.  The distinguished order unit of $K_0\vep$ is the
constant polynomial 1.\label{dgproposition}
\end{theo}

\begin{proof}
We will construct an order isomorphism from $K_0\vep$ into $G$. The
transposes of the incidence matrices will be used for typographical
reasons in the computation in order to make the computations on row
vectors. For $p(x)=a_0+\dots a_dx^d$, the transpose of the $k$'th incidence
matrices associated to $\vep$, $\phi_k$ will be a $((k-1)d+1) \times (kd+1)$ matrix
with \[(\phi_k)_{ij}=\left\{\begin{array}{ll} a_{(j-i)}& \text{ if
}0\leq
j-i \leq d\\
0 & \text{otherwise}\end{array}\right.\]

For $l\leq m$, define $\phi_{lm}:\Z^{d(l-1)+1}\to\Z^{dm+1}$ by
$\phi_l\phi_{l+1}\dots\phi_{m}$.  We will identify $\Z^i$ with the
additive group of polynomials of degree at most $i-1$, $\Z_{i-1}[x]$
in the following manner.  For $v=[v_0\ v_1\dots v_{i-1}]\in \Z^i$,
define $v(x)\in \Z_{i-1}[x]$ by $v(x)=\sum_{j=0}^{i-1}v_jx^j$.
Now if $v\in \Z^{dm+1}$, we have $(v\phi_m)(x)=v(x)p(x)$. Under the
above correspondence, $\phi_l$ becomes multiplication by $p(x)$ for
all $l$, and $\phi_{lm}$ becomes multiplication by $(p(x))^{m-l}$.

Define $\rho_m:\Z_{md}[x]\to G$ by $\rho_m(r(x))=\ds
\frac{r(x)}{(p(x))^m}$. In order to satisfy the hypothesis of the
universal mapping property of direct limits, it needs to be shown
that for $l\leq m$, $\rho_l=\rho_m\circ\phi_{lm}$:

\[\begin{array}{ccl}
  \rho_m\circ\phi_{lm}(r(x))
  &=&\rho_m(r(x)(p(x))^{m-l})\\
  &=&\ds\frac{r(x)(p(x))^{m-l}}{(p(x))^m}\\
  &=&\ds\frac{r(x)}{(p(x))^l}\\
  &=&\rho_l(r(x)).
  \end{array}\]

  Hence the hypothesis for the universal mapping property of direct limits is satisfied,
and the $\rho_l$ are constant on equivalence classes.  It follows
that there is a unique homomorphism $\rho:K_0\vep\to \G$, which can
be defined on an equivalence class by taking any representative in
$\Z^{di+1}$ and applying $\rho_i$ to it.  This is well defined
because $\rho_i$ is constant on equivalence classes, and there is
only one element of each equivalence class in each
  $\Z^{id+1}$.
  We claim that $\rho$ is an isomorphism.

  First we show that $\rho$ is a homomorphism.  For without loss of generality, assume $l\leq
  m$, $r(x)\in \Z_{md}[x]$, and $s(x)\in \Z_{ld}[x]$.  Then
  \[\begin{array}{ccl}
  \rho\left(\overline{r(x)}+\overline{s(x)}\right) &=&
  \rho\left(\overline{r(x)+(p(x))^{m-l}s(x)}\right)\\
  &=& \ds \frac{r(x)+(p(x))^{m-l}s(x)}{(p(x))^m}\\
    &=&\ds \frac{r(x)}{(p(x))^m} + \frac{s(x)}{(p(x))^l}\\
  &=& \rho(\overline{r(x)})+\rho(\overline{s(x)}).
  \end{array}\]

Now we show that $\rho$ is onto.  Given $\ds \frac{r(x)}{(p(x))^m}\in
\G$, then for $r(x)\in \Z_{md}[x]$,
$\rho(\overline{r(x)})=\ds \frac{r(x)}{(p(x))^m}$.\\

Lastly we show that $\rho$ is injective. If $r(x)\in\Z_{md}[x]$,
$\rho(\overline{r(x)})=0$, then
$$ \frac{r(x)}{(p(x))^m}=0\text{ therefore } r(x)=0\text{ and }
\overline{r(x)}=\overline{0}.$$

Hence $\rho$ is an isomorphism, and $\G$ is isomorphic to $K_0\vep$.
In addition, $\G$ is order isomorphic to $K_0\vep$ because
$$(\G)_+=\left\{\frac{r(x)}{(p(x))^m}|\ r(x)(p(x))^l\text{ has all positive
coefficients for some } l\geq 0\right\}$$ is exactly the image of
the positive set of $\ds\lim_{\rightarrow}\Z^{dk+1}$ under $\rho$.
Finally, the image of 1 is $\ds \frac{1}{(p(x))^0}=1$.
\end{proof}

\section{Some Ergodic Adic Invariant Measures}
  Determining the ergodic adic invariant measures for systems in
  $\Sl$ depends heavily on the particular system.  For instance the
  Euler adic has a unique fully supported (every cylinder set is given positive measure) invariant measure,
  \cite{BKPS, BP} while the reverse Euler, \cite{FP}, and the polynomial systems have one-parameter
  families of ergodic adic invariant measures.  We will devote the
  remainder of this section to showing that each of the adics given by
  a positive integer polynomial has a one-parameter family of ergodic adic invariant measure.

For a cylinder set $C$ in $\X$, and any path $\gamma\in \X$,
$\dim(C,(n,k_n(\gamma)))$ is the number of paths from the terminal
vertex $(m,l)$ of $C$, to the vertex $(n,k_n(\gamma))$.  Define a
function $\co: \Z\times\Z\to\Z$ by $\co(n,k)=$ the coefficient of
$x^{k}$ in the polynomial $(p(x))^n$. Because of the self-similarity
of this class of Bratteli diagrams, if $C$ terminates at vertex
$(m,l)$,
\begin{align*}\dim(C,(n,k_n(\gamma)))&=\co(n-m,k_n(\gamma)-l)\text{ and }\\
\dim(n,k)&=\co(n,k).\end{align*}

\begin{figure}[h]
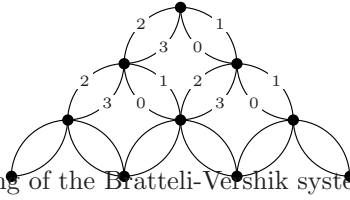

\begin{center}
\begin{graph}(5,1.5)(0,1.5)
\unitlength=.75\unitlength
\roundnode{30}(1,1)\roundnode{31}(3,1)\roundnode{32}(5,1)
\roundnode{33}(7,1)\roundnode{20}(2,2)\roundnode{21}(4,2)\roundnode{22}(6,2)
\roundnode{10}(3,3)\roundnode{11}(5,3)\roundnode{00}(4,4)
\bow{00}{11}{.2}\bow{00}{11}{-.2}\bow{00}{10}{.2}
\bow{00}{10}{-.2}\bow{10}{20}{.2}\bow{10}{20}{-.2}\bow{10}{21}{.2}\bow{10}{21}{-.2}
\bow{11}{21}{.2}\bow{11}{21}{-.2}\bow{11}{22}{.2}\bow{11}{22}{-.2}
\bow{20}{30}{.2}\bow{20}{30}{-.2}\bow{20}{31}{.2}\bow{20}{31}{-.2}
\bow{21}{31}{.2}\bow{21}{31}{-.2}\bow{21}{32}{.2}\bow{21}{32}{-.2}
\bow{22}{32}{.2}\bow{22}{32}{-.2}\bow{22}{33}{.2}\bow{22}{33}{-.2}
\bowtext{00}{10}{-.2}{\tiny{2}}\bowtext{00}{10}{.2}{\tiny{3}}
\bowtext{00}{11}{-.2}{\tiny{0}}\bowtext{00}{11}{.2}{\tiny{1}}
\bowtext{10}{20}{-.2}{\tiny{2}}\bowtext{10}{20}{.2}{\tiny{3}}
\bowtext{10}{21}{-.2}{\tiny{0}}\bowtext{10}{21}{.2}{\tiny{1}}
\bowtext{11}{21}{-.2}{\tiny{2}}\bowtext{11}{21}{.2}{\tiny{3}}
\bowtext{11}{22}{-.2}{\tiny{0}}\bowtext{11}{22}{.2}{\tiny{1}}
\end{graph}
\end{center}
\caption{Labeling of the Bratteli-Vershik system
$\VE_{2+2x}$}\label{labelpic}
\end{figure}

For  $n=0,1,\dots$ and $d\leq k\leq d(n-1)$, the number of edges
into vertex $(n,k)$ is exactly $a_0+a_1+\dots +a_d$.  In addition,
for every vertex $(n,k)$ the number of edges leaving $(n,k)$ is
exactly $a_0+\dots+a_d$. Because of this it is convenient to use an
alphabet to label the edges of paths in $\X$.  The alphabet
associated to $\X$ will be $A=\{0,1,\dots, a_0+a_1+\dots +a_d-1\}$.
If an edge $e$ is the $j$'th edge between vertex $(n,k)$ and $(n+1,k+i)$ label it
$$\left(\sum_{m=0}^{d-i-1}\right)+(j-1);$$  see Figure \ref{labelpic}.  By labeling in this manner,
the lexicographic ordering on comparable edges is consistent with
the edge ordering given for the general family $\Sl$.  Then any path
in $\X$ is uniquely determined by the labeling of its edges, and
because of this, we use both $X$ and a one sided infinite sequence
in $A^{\N}$ to denote the infinite edge paths on a Bratteli diagram.
For ease of notation we will refer to a path by the infinite
labeling of its edges, and when the context is clear, we will refer
to an edge by its label.

Consider some cylinder set $C=[c_0c_1...c_{n-1}] \in \X$ and any
$\T-$invariant Borel probability measure $\mu$ on $\X$.  Define the
\emph{weight} $w_{c_0}$ on the edge $c_0$ to be $\mu([c_0])$. For
$n>0$ and $\mu ([c_0c_1\dots c_{n-1}])=0$ define the weight
$w_{c_n}$ on $c_n$ to be 0.  For $n>0$ and $\mu ([c_0\dots
c_{n-1}])>0$ define $w_{c_n}$ on $c_n$ to be
$\mu([c_0...c_n])/\mu([c_0...c_{n-1}])$.  Then
$\mu([c_0...c_n])=w_{c_0}...w_{c_n}$.

These weights are well defined because as we will see in Lemma
\ref{claim1}, all cylinders with the same terminal vertex have the
same measure.

\begin{rem}In this section we discuss
measures, $\T$-invariant Borel probability measures for which edges
with the same label have the same weight. Then for the probability
space $(\X,\mathcal{B},\mu)$ there are at most $a_0+a_1+\dots +a_d$
different weights.  For each $j\in A$ we will denote by $w_j$ the
weight associated to each edge labeled $j$.  Since
$(\X,\mathcal{B},\mu)$ is a probability space,
$\sum_{i=0}^{a_0+\dots+a_d-1}w_i=1$.  In view of the labeling of
paths by the alphabet $A$, these measures are Bernoulli and we
denote such a measure by $\mathcal{B}(w_{a_0+\dots +a_d-1},\dots ,
w_0)$. \end{rem}

\begin{figure}[h]
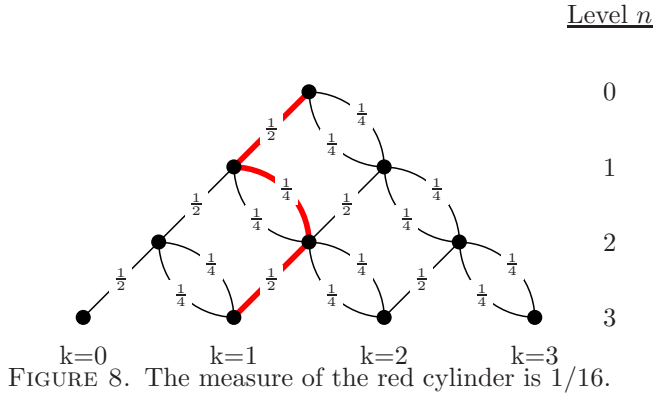

\begin{center}
\begin{graph}(4,4)(0,1)
\unitlength=\unitlength
\roundnode{00}(2,4)\roundnode{10}(1,3)\roundnode{11}(3,3)
\roundnode{20}(0,2)\roundnode{21}(2,2)\roundnode{22}(4,2)
\freetext(6,5){\underline{Level $n$}} \freetext(6,4){0}
\freetext(6,3){1}\freetext(6,2){2}\freetext(6,1){3}
\roundnode{30}(-1,1)
\roundnode{31}(1,1)\roundnode{32}(3,1)\roundnode{33}(5,1)
\edge{00}{10}[\graphlinecolour(1,0,0)\graphlinewidth{.075}]\bow{00}{11}{-.2}\bow{00}{11}{.2}
\edge{10}{20}\bow{10}{21}{-.2}
\bow{10}{21}{.2}[\graphlinecolour(1,0,0)\graphlinewidth{.075}]
\edge{11}{21}\bow{11}{22}{-.2}\bow{11}{22}{.2}
\edge{20}{30}\bow{20}{31}{-.2}\bow{20}{31}{.2}
\edge{21}{31}[\graphlinecolour(1,0,0)\graphlinewidth{.075}]\bow{21}{32}{-.2}\bow{21}{32}{.2}
\edge{22}{32}\bow{22}{33}{-.2}\bow{22}{33}{.2}
\nodetext{30}(0,-.5){k=0}\nodetext{31}(0,-.5){k=1}\nodetext{32}(0,-.5){k=2}\nodetext{33}(0,-.5){k=3}
\edgetext{00}{10}{\tiny$\frac{1}{2}$}\bowtext{00}{11}{-.2}{\tiny$\frac{1}{4}$}\bowtext{00}{11}{.2}{\tiny$\frac{1}{4}$}
\edgetext{10}{20}{\tiny$\frac{1}{2}$}\bowtext{10}{21}{-.2}{\tiny$\frac{1}{4}$}\bowtext{10}{21}{.2}{\tiny$\frac{1}{4}$}
\edgetext{11}{21}{\tiny$\frac{1}{2}$}\bowtext{11}{22}{-.2}{\tiny$\frac{1}{4}$}\bowtext{11}{22}{.2}{\tiny$\frac{1}{4}$}
\edgetext{20}{30}{\tiny$\frac{1}{2}$}\bowtext{20}{31}{-.2}{\tiny$\frac{1}{4}$}\bowtext{20}{31}{.2}{\tiny$\frac{1}{4}$}
\edgetext{21}{31}{\tiny$\frac{1}{2}$}\bowtext{21}{32}{-.2}{\tiny$\frac{1}{4}$}\bowtext{21}{32}{.2}{\tiny$\frac{1}{4}$}
\edgetext{22}{32}{\tiny$\frac{1}{2}$}\bowtext{22}{33}{-.2}{\tiny$\frac{1}{4}$}\bowtext{22}{33}{.2}{\tiny$\frac{1}{4}$}
\end{graph}
\end{center}
\caption[An example of a $\T$-invariant measure.]{The measure of the
red cylinder is $1/16$.}\label{berexam}
\end{figure}

In \cite{Xavier, Mela}, M\'{e}la showed that when all the coefficients of
$p(x)$ are 1,  the invariant ergodic probability measures for
$T_{p(x)}$ are the Bernoulli measure $\B(0,...,0,1)$ and the
one-parameter family $\B(q,t_q,\ds
\frac{t_q^2}{q},\frac{t_q^3}{q^2},..,\frac{t_q^n}{q^{n-1}})$, where
$t_q$ is the unique solution in [0,1] to the equation
$$q^{n}-q^{n-1}+q^{n-1}t+q^{n-2}t^2+...+qt^{n-1}+t^n =0.$$

 Using similar techniques we have extended this result to the following:

\begin{theo}  Let $p(x)=a_0+\dots a_dx^d$ and let $(\X,\T)$ be the
Bratteli-Vershik system in $\Slp$ determined by $p(x)$. If $q\in
\left(0,\ds \frac{1}{a_0}\right)$, and $t_q$ is the unique solution in $[0,1]$
to the equation
\begin{equation}
a_0q^d+a_1q^{d-1}t+...+a_dt^d-q^{d-1}=0,
\end{equation}
then the invariant, fully supported, ergodic probability measures
for the adic transformation $T_{p(x)}$ are the one-parameter family
of Bernoulli measures

$$\B\left(\underbrace{q,...,q}_{a_0\
 times},\underbrace{t_q,...,t_q}_{a_1\ times},\underbrace{\ds
\frac{t_q^2}{q},...,\frac{t_q^2}{q}}_{a_2\ times},\dots
,\underbrace{\frac{t_q^n}{q^{n-1}},..., \frac{t_q^n}{q^{n-1}}}_{a_n\
times}\right).$$\label{iemab}
\end{theo}

\begin{prop}   Let $p(x)=a_0+\dots a_dx^d$ and $(\X,\T)$ be the
Bratteli-Vershik system determined by $p(x)$. The only $\T$
invariant, ergodic probability measures that are not fully supported
are the Bernoulli measures
$$\B\left(\ds\underbrace{\frac{1}{a_0},\dots,\frac{1}{a_0}}_{a_0\
times},0,\dots,0\right)\text{ and }\B\left(0,
\dots,0,\ds\underbrace{\frac{1}{a_n},\dots,\frac{1}{a_n}}_{a_n\
times}\right).$$\label{nfsthm}
\end{prop}

The proofs of Theorem \ref{iemab} and Proposition \ref{nfsthm} use
many other results and definitions, which are presented below.
Proposition \ref{eab} shows that every invariant fully supported
ergodic probability measure for $(\X,\T)$ must be Bernoulli.
Proposition \ref{bre} says that the Bernoulli measures that are
$\T$-invariant are in fact ergodic. Proposition \ref{bai} shows
which Bernoulli measures are $\T$-invariant. This will prove Theorem
\ref{iemab}. We will then conclude with the proof of Proposition
\ref{nfsthm}.

\begin{lemma}
Any non-atomic measure on $\X$ is $\T$-invariant if and only if all
cylinders with the same terminal vertex have the same measure.
\label{claim1}
\end{lemma}
\begin{proof}
($\Rightarrow$) Let $(n,k)\neq (0,0)$ be a vertex of
$(\mathcal{V,E})_{p(x)}$.  Consider the maximal path from (0,0) to
$(n,k)$ and the cylinder set, $C_{\max}$, defined by this path.
There exists an $i$ such that $\T^{-i}(C_{\max})$ is the minimal
cylinder $Y_n(k,0)$ determined by the minimal path from $(0,0)$ to
$(n,k)$. Since the measure is $\T$-invariant, the elements of the
set $\{\T^{-j}(C_{\max})\}_{j=0}^i$ all have the same measure.
Because all the cylinders with terminal vertex $(n,k)$ are contained
in the above set, all cylinders with
terminal vertex $(n,k)$ have the same measure.

($\Leftarrow$)  It is enough to show that for each cylinder set $C$,
$\T^{-1}C$ has the same measure as $C$.  Let $C$ be any cylinder
set, with terminal vertex $(n,k)$.   Suppose that $C$ is not minimal.  Since $C$ is not minimal,
$\T^{-1}(C)$ has terminal vertex $(n,k)$, and hence the same measure
as $C$.

If $C$ is minimal, it can be decomposed into a disjoint union of minimal cylinders and at most a countable
number of infinite paths.  since the measure is non-atomic, the measure of $C$ equals the measure of $T^{-1}C$.
\end{proof}

The following is a lemma of Vershik that is proved using the ergodic
theorem and the indicator function for the cylinder set $C$.
\begin{lemma}[Vershik \cite{Vershik1974, VerKer1981}]  If $\mu$ is an invariant non-atomic ergodic
probability measure for the adic transformation $T_{p(x)}$, then for
every
cylinder set $C$, \\
\centerline{$\ds \mu(C)= \lim_{n\to
\infty}\frac{\dim(C,(n,k_n(\gamma))}{\dim(n,k_n(\gamma))}$ for
$\mu$-a.e. $\gamma\in X$}\label{lemma1}
\end{lemma}

The following lemma shows that the invariant and ergodic probability
measures for $\T$ are also invariant for the one-sided shift
$\sigma$ on $A^{\N}$.

\begin{lemma} For each $\gamma\in \X$ and $j\in A$,
let $\sigma_j\gamma=j\gamma_0\gamma_1...$.  If $\mu$ is invariant
and ergodic for $\T$, then for any cylinder set $C$,
\centerline{$\mu(C)=\mu(\sigma_0C)+\mu(\sigma_1C)+...+\mu(\sigma_{
a_d+...+a_0-1}C)=\mu(\sigma^{-1}C).$}\label{lemma2}
\end{lemma}
\begin{proof}
Define a function $g:A\to \{0,1,...,d\}$ such that if the letter $j\in A$ is the label of an edge which
connects vertex $(n,k)$ to $(n,k+i)$, then $g(j)=i$.  Then for any cylinder set $C=[c_0...c_{n-1}]$ which
terminates at vertex $(n,k)$, we have $\ds \sum_{i=0}^{n-1}
g(c_i)=k$.

Let $C$ be a cylinder set with terminal vertex $(m,l)$.  For almost
every $\gamma$ in $X$ and each $j\in A$ we have
$$\mu(\sigma_jC)= \lim_{n\to
\infty}\frac{\dim(\sigma_jC,(n,k_n(\gamma)))}{\dim(n,k_n(\gamma))}.$$
Denote by $C^j$ the cylinder set extended by $j$.  The terminal vertex of $C^j$ is $(m+1,l+g(j)),$ and since the first
$m+1$ edges of $\sigma_jC$ are a permutation of those of $C^j$, the
terminal vertex of $\sigma_jC$ is also $(m+1,l+g(j))$.  Hence for
all $n>m$,
$\dim(C^j,(n,k_n(\gamma)))=\dim(\sigma_jC,(n,k_n(\gamma)))$.

The set of finite paths starting from $(m,l)$ and ending at
$(n,k_n(\gamma))$ can be divided into $a_d+...+a_0$ groups,
according to whether the edge is labeled $0, 1, \dots \ds
\left(\sum_{i=0}^{d}a_i\right)-1$. Then we have, $\dim(C,(n,k_n(\gamma)))$
\[\begin{array}{rll}
 & = \dim(C^0,(n,k_n(\gamma))) +...+
\dim(C^{(\sum_{i=0}^{d}a_i)-1},(n,k_n(\gamma)))\\
&=\dim(\sigma_0C,(n,k_n(\gamma)))+...+\dim(\sigma_{(\sum_{i=0}^{d}a_i)-1}C,(n,k_n(\gamma))).\end{array}\]
Therefore, $\ds \frac{\dim(C,(n,k_n(\gamma)))}{\dim(n,k_n(\gamma))}$\\$$=\ds
\frac{\dim(\sigma_0C,(n,k_n(\gamma)))}{\dim(n,k_n(\gamma))}+...+\ds
\frac{\dim(\sigma_{(\sum_{i=0}^{d}a_i)-1}C,(n,k_n(\gamma)))}{\dim(n,k_n(\gamma))}.$$
Taking limits as $n\to \infty$, $\mu(C)=\mu(\sigma_0C)+...+\mu(\sigma_{a_d+...+a_0-1}C).$
\end{proof}
\begin{lemma}For $j_0,j_1\in A$,
$\dim(C^{j_0},(n,k_n(\gamma)))=\dim(C^{j_1j_0},(n,k_{n+1}(\sigma_{j_1}\gamma))).$
\end{lemma}
\begin{proof}
Assume that $C$ terminates at vertex $(m,k)$.  Then $C^{j_0}$
terminates at $(m+1,l+g(j_0))$, where $g(j)$ is as in the proof of
Lemma \ref{lemma2}. Hence, $\dim(C^{j_0},(n,k_n(\gamma)))=\co
(n-(m+1),k_n(\gamma)-(l+g(j_0))$. Also,
$k_{n+1}(\sigma_{j_1}\gamma)=k_n(\gamma)+g(j_1)$, and $C^{j_1j_0}$
terminates at vertex $(m+2,l+g(j_1)+g(j_0))$.  Hence
$$\dim(C^{j_1j_0},(n,k_{n+1}(\sigma_{j_1}\gamma)))=\co (n+1-(m+2),k_n(\gamma)+g(j_1)-(l+g(j_1)+g(j_0)))$$
$$=\co(n-(m+1),k_n(\gamma)-(l+g(j_0)))$$
$$=\dim(C^{j_0},(n,k_n(\gamma))).$$
\end{proof}

\begin{prop} Every $\T$-invariant fully supported ergodic probability measure
for $(\X,\T)$ is Bernoulli.\label{eab}
\end{prop}

\begin{proof}  Let $\mu$ be a $\T$-invariant fully supported
ergodic probability mesure for $(\X,\T)$.  To prove that $\mu$ is a
Bernoulli measure, it is enough to show that for each $i\in A$ there
exists a number $p_i$ such that for every cylinder set $C$,
$\ds\frac{\mu(C^i)}{\mu(C)}=p_i$. Now for any $\gamma\in \X$,

$$\frac{\dim(C^i,(n,k_n(\gamma)))}{\dim(C,(n,k_n(\gamma)))}=\ds
\frac{\dim(C^{ji},(n,k_{n+1}(\sigma_j\gamma)))}{\dim(C^j,(n,k_{n+1}(\sigma_j\gamma)))}.$$

By Lemma \ref{lemma1}, there exists a set $E$ of full measure such
that for all $\gamma\in E$ and all $i\in A$
$$\ds\frac{\mu(C^i)}{\mu(C)}=\lim_{n\to
\infty}\frac{\dim(C^i,(n,k_n(\gamma)))}{\dim(C,(n,k_n(\gamma)))}.$$

If there is $j\in A$ such that $E\cap \sigma_jE=\emptyset$, then
$\mu(E\cap \sigma_jE)=0$, and, since $E$ has full measure,
$\mu(\sigma_jE)=0$. Denote by $[r]$ the cylinder set $\{\gamma\in
X:\gamma_0=r\}$.

Then by Lemma \ref{lemma2}:
$$1=\mu(E)=\sum_{r\neq j}\mu(\sigma_rE)\leq
\sum_{r\neq j}\mu[r]\leq 1\text{ which implies } \mu([j])=0,$$
contradicting our earlier assumption that $\mu$ has full support.
Hence there exists $\gamma\in E\cap\sigma_jE$.  Let $\xi$ be the
path in $E$ such that $\sigma_j\xi=\gamma$; then
$$\lim_{n\to
\infty}\frac{\dim(C^i,(n,k_n(\xi)))}{\dim(C,(n,k_{n}(\xi)))}=\lim_{n\to\infty}
\frac{\dim(C^{ji},(n,k_{n+1}(\sigma_j\xi)))}{\dim(C^j,(n,k_{n+1}(\sigma_j\xi)))},$$
$$\text{showing that } \frac{\mu(C^i)}{\mu(C)}=\frac{\mu(C^{ji})}{\mu(C^j)}.$$

Then for any cylinder set $C=[c_0c_1...c_{m-1}]$ we have:
$$\frac{\mu(C^i)}{\mu(C)}=\frac{\mu([c_0...c_{m-2}]^{c_{m-1}i})}{\mu([c_0...c_{m-2}]^{c_{m-1}})}
=\frac{\mu([c_0...c_{m-2}]^{i})}{\mu([c_0...c_{m-2}])}=...=\frac{\mu([c_0]^i)}{\mu([c_0])}.$$

Also, for all $j,k,l\in A$ we have:\\
$\mu([jli])=\mu([lji])$, since $[jli]$ and $[lji]$ have the same
terminal vertex.  Then
$$\frac{\mu([j]^{li})}{\mu([j]^l)}=\frac{\mu([l]^{ji})}{\mu([l]^j)}\text{ so that }
\frac{\mu([j]^{i})}{\mu([j])}=\frac{\mu([l]^{i})}{\mu([l])}.$$ This
shows that $\ds \frac{\mu(C^i)}{\mu(C)}$ is independent of $C$, and
hence equal to $\mu([i])$.  Therefore $\mu$ is a Bernoulli.
\end{proof}

\begin{prop} The $\T$-invariant Bernoulli measures on $\X$
are ergodic.\label{bre}
\end{prop}
\begin{proof}
Define the random variable $Z_i$ on $\X$ by letting $Z_i(\gamma)$ be
the label on the $i-1$'th edge of $\gamma$.  Since the probability measure is
Bernoulli, the $Z_i$ are independent and identically distributed. If
$B$ is a set that depends symmetrically on $Z_1,\dots,Z_n$, then
$\gamma\in B$ implies that $\{\xi\in X|\xi_0\xi_1\dots \xi_{n-1}
\text{ is a permutation of }\gamma_0\gamma_1\dots \gamma_{n-1}\text{
and for }m\geq n,\ \xi_m=\gamma_m\}$ is also in $B$. If $\mathcal{S}_n$ is $\sigma$-algebra set generated by
such $B$, the Hewitt-Savage theorem implies that $\ds
\mathcal{S}=\cap_{n=1}^\infty \mathcal{S}_n$ is trivial.

Let $\mathcal{T}_n$ be the $\sigma$-algebra generated by sets $B'$
such that if $\gamma\in B'$, then $\{\xi\in X|\text{ for }m\geq n,\
\xi_m=\gamma_m\}$ is also in $B'$.  Then for
each generator $B'$ of $\mathcal{T}_n$, there are a finite number of
generators $B_i$ of $\mathcal{S}_n$ such that
$\cup_{i=1}^{m}B_i=B'$. Hence $B'\subset \mathcal{S}_n$.  Then
$\mathcal{T}_n\subset \mathcal{S}_n$, and
$\cap_{i=1}^{\infty}\mathcal{T}_n=\mathcal{T}\subset\mathcal{S}$.
Since $\mathcal{S}$ is trivial, so is $\mathcal{T}$.  But
$\mathcal{T}$ is the $\sigma$-algebra of $\T$-invariant sets.
Therefore the invariant Bernoulli measures for $\T$ are ergodic.
\end{proof}

It remains to determine which Bernoulli measures are invariant.

\begin{prop}  The Bernoulli measures invariant for the adic
transformation $\T$ are the fully supported ones described in
Theorem \ref{iemab}, along with
$$\B\left(\ds\underbrace{\frac{1}{a_0},\dots,\frac{1}{a_0}}_{a_0\
times},0,\dots,0\right)\text{ and }\B\left(0,
\dots,0,\ds\underbrace{\frac{1}{a_d},\dots,\frac{1}{a_d}}_{a_d\
times}\right).$$\label{bai}
\end{prop}
\begin{proof}
  Recall that any edge label $j$ has weight $w(j)$.
 Recall the definition of $g:A\to \{0,1,\dots,d\}$
  as given in the proof of Lemma \ref{lemma2}.  By
  Lemma \ref{claim1}, $g(j_1)=g(j_2)$ implies $w(j_1)=w(j_2)$.  For
  $0\leq t\leq d$, define $p_t=w(j)$ whenever $g(j)=t$.
Then \begin{equation} a_0p_0+...+a_dp_d=1. \label{1}\end{equation}
For $s\in \N$ and $i_k,j_k\in \{0,1,\dots,d\}$ for all
$k=0,1,\dots,s$, Lemma \ref{claim1} implies that a Bernoulli measure
is $\T$ invariant if and only if whenever
\begin{equation}
\sum_{k=0}^{s}i_k=\sum_{k=0}^{s}j_k, \text{ we have
}\prod_{k=0}^{s}p_{i_k}=\prod_{k=0}^{s}p_{j_k}.\label{product}\end{equation}

Assume for now that $p_0,p_1>0$. Claim:  Equation \ref{product} is
satisfied if and only if $p_0p_j=p_1p_{j-1}$ for $1\leq j\leq d$.

Clearly Equation \ref{product} implies $p_0p_j=p_1p_{j-1}$. It
remains to be shown that $p_0p_j=p_1p_{j-1}$ implies Equation
\ref{product}.

For $1\leq j\leq d$ we will assume
\begin{equation}p_0p_j=p_1p_{j-1}.\label{claimequ}\end{equation}
We will use induction to prove to prove our claim.  The hypothesis
is that for \\$i_0,i_1\dots, i_{s-1}, j_0,j_1,\dots,j_{s-1}$ in
$\{0,1,\dots,d\}$, whenever
\begin{equation}i_0\dots +i_{s-1}  = j_0+\dots
j_{s-1},\text{ we have }
\prod_{k=0}^{s-1}p_{i_k}=\prod_{k=0}^{s-1}p_{j_k}.\end{equation} We
will show that for $i_0,i_1\dots, i_{s}, j_0,j_1,\dots,j_{s}$ in
$\{0,1,\dots,d\}$, whenever
\begin{equation*}i_0\dots +i_{s}  = j_0+\dots
j_{s},\text{ we have }
\prod_{k=0}^{s}p_{i_k}=\prod_{k=0}^{s}p_{j_k}.\end{equation*}

We now show the base case. For $1\leq i\leq d$ and $0\leq k\leq
d-1$, Equation \ref{claimequ} implies
$$p_i=\frac{p_1}{p_0}p_{i-1}\text{ and
}p_k=\frac{p_0}{p_1}p_{k+1}.$$ Hence
$$p_ip_k=\frac{p_1}{p_0}p_{i-1}\frac{p_0}{p_1}p_{k+1}=p_{i-1}p_{k+1}.$$
For $i,k,l,m\in \{0,1,\dots,d\}$ we then have that whenever
$i+k=l+m$, $p_ip_k=p_lp_m$, hence we have shown the base case.

Now consider $i_0,i_1,\dots,i_s,j_0,j_1,\dots,j_s$ in
$\{0,1,\dots,d\}$ such that
$$i_0+\dots+i_s=j_0+\dots+j_s.$$
Then $i_0+\dots+i_s-i_s-j_s=j_0+\dots+j_s-j_s-i_s$, hence
$$i_0+\dots+i_{s-1}-j_s=j_0+\dots+j_{s-1}-i_s.$$  There also exist
$l_0,l_1,\dots,l_{s-2}$ in $\{0,1,\dots,d\}$ such that
$$l_0+\dots+l_{s-2}=i_0+\dots+i_{s-1}-j_s=j_0+\dots+j_{s-1}-i_s.$$
Adding $j_s$ to both sides, we see that
$l_0+\dots+l_{s-2}+j_s=i_0+\dots+i_{s-1}$, hence the induction
hypothesis implies
\begin{equation}p_{j_s}\prod_{k=0}^{s-2}p_{l_k}=\prod_{k=0}^{s-1}p_{i_k}.\label{eq1}\end{equation}
Likewise, $l_0+\dots+l_{s-2}+i_s=j_0+\dots+j_{s-1}$, and the
induction hypothesis implies
\begin{equation}p_{i_s}\prod_{k=0}^{s-2}p_{l_k}=\prod_{k=0}^{s-1}p_{j_k}.\label{eq2}\end{equation}
Combining Equation \ref{eq1} and Equation \ref{eq2} we see that
\begin{equation*}
\prod_{k=0}^{s}p_{i_k}=p_{i_s}p_{j_s}\prod_{k=0}^{s-2}p_{l_k}
=p_{j_s}p_{i_s}\prod_{k=0}^{s-2}p_{l_k}=\prod_{k=0}^{s}p_{j_k}.
\end{equation*}
Therefore we have proved the claim and the Bernoulli measures are
$\T$ invariant if and only if for $1\leq j \leq d$,
$$p_0p_j=p_1p_{j-1}.$$

For simplicity of notation define $p_0=q$ and $p_1=t$.  For $1\leq
j\leq d$, $p_j=\dfrac{t}{q}p_{j-1}$.  Hence every $p_j$ can be
defined inductively by $t$ and $q$.  In particular, for $1\leq j\leq
d$,
$$p_j=\frac{t^j}{q^{j-1}}.$$
By Equation \ref{1}
\begin{equation*}
a_0q+a_1t+a_2\frac{t^2}{q}+\dots+a_d\frac{t^d}{q^{d-1}}=1.
\end{equation*}

Multiplying through by $q^{d-1}$ and simplifying, we see
that\begin{equation}a_0q^d+a_1q^{d-1}t+...+a_dt^d-q^{d-1}=0.\label{poly}\end{equation}

To conclude that $p_1,...,p_d$ are completely determined by the
choice of $p_0=q$, it remains only to show that for each $q\in
(0,1/a_0)$, Equation \ref{poly} has a unique solution in [0,1].

Consider $m(t)=a_0q^d+a_1q^{d-1}t+...+a_dt^d-q^{d-1}$.  Then
$m(0)=a_0q^d-q^{d-1}=q^{d-1}(a_0q-1)\leq 0$, since $a_0q\leq 1$.
Also, $m(1)=a_0q^d+a_1q^{d-1}+...+a_d-q^{d-1}>0$, since $a_1\geq 1$
implies $a_1q^{d-1}-q^{d-1}\geq 0$.  By the intermediate value
theorem, there exists a root in $[0,1]$. Now,
$m'(t)=a_1q^{d-1}+...+da_dt^{d-1}$ is strictly greater than 0 on
$[0,1]$, so that $m(t)$ is strictly increasing on $[0,1]$; therefore
there is a unique solution $t_q$ to $m(t)=0$ in the interval
$[0,1]$.

If $p_0=1/a_0$, $a_0p_0=1$ and all other $p_i=0$, hence the
$\T$-invariant measure $\mu$ is supported on the paths for which
$k_n(\gamma)=0$ for all $n\geq 0$.  Finally, if $i\leq d-2$ and
$p_i=0$, then $i+2\leq d$, and $p_ip_{i+2}=p_{i+1}p_{i+1}=0$. Hence
$p_{i+1}=0$. If $p_0=0$ then $p_i=0$ for all $0\leq i<d$. This
implies that the only nonzero probability is $p_d$ and $a_dp_d=1$,
hence $p_d=1/a_d$ and the $\T$-invariant measure $\mu$ is supported
on the set of paths for which $k_n(\gamma)=dn-1$ for all $n\geq 0$.

\end{proof}

We now have enough tools to prove Theorems \ref{iemab} and \ref{nfsthm}.

\begin{proof}[Proof of Theorem \ref{iemab}]  Direct from Proposition \ref{bai}, Proposition \ref{bre}, and Proposition
\ref{eab}.
\end{proof}

\begin{proof}[Proof of Proposition \ref{nfsthm}]
According to Proposition \ref{bai} and Proposition \ref{bre}, the
Bernoulli measures
$$\B\left(\ds\underbrace{\frac{1}{a_0},\dots,\frac{1}{a_0}}_{a_0\
times},0,\dots,0\right)\text{ and }\B\left(0,
\dots,0,\ds\underbrace{\frac{1}{a_n},\dots,\frac{1}{a_n}}_{a_n\
times}\right)$$ are ergodic and invariant for $(\X,\T)$.  It remains
to show that these are the only invariant, ergodic probability
measures that are not fully supported on $\X$.

 By Proposition
  \ref{prop311}, the only proper closed invariant sets are those for
  which the tails are eventually diagonal.,

Define $A_k$ to be the closed invariant set $$A_k=\{\gamma\in
X|\text{ either }k_n(\gamma)\leq k \text{ for all } n \text{ or
}dn-k_n(\gamma)\leq k\text{ for all }n\}.$$ Define $C$ to be the maximal cylinder
from the root vertex to vertex $(l,k)$, where $l $ is the first
level for which $dl-k> k$. In other words, $C$ is the maximal
cylinder of shortest length such that for each $\gamma\in C$ and
$n\geq l$, $k_n(\gamma)=k$.

If $\mu$ is a $\T$-invariant ergodic probability measure that is not
fully supported on $\X$, we will compute $\mu(C)$ using Lemma
\ref{lemma1}. If there is a set of positive measure $B_k\subset A_k$
for which every $\gamma\in B_k$ and $n\geq l$ has $k_n(\gamma)\neq
k$, then for all $n\geq l$,  $\dim(C,(n,k_n(\gamma)))= 0$, hence
$\mu(C)=0$. Therefore we will assume that almost every $\gamma$ in
$A_k$ has $k_n(\gamma)\equiv k$ for sufficiently large $n$.  This
implies that the measure $\mu$ is supported on the paths for which
$k_n(\gamma)\leq k$.

Then
$$\mu(C)=\lim_{n\to\infty}\frac{\dim(C,(n,k_n(\gamma)))}{\dim(n,k_n(\gamma))}.$$
It is clear that $\dim(C,(n,k_n(\gamma)))= a_0^{n-l}$.  We will now
find a lower bound for $\dim(n,k_n(\gamma))$.  Let $m\in \Z_+$ such
that $k-md>0$ and $k-(m+1)d\leq 0$.  Then $dm+i=k$ for some $k\in
\{1,2,\dots, d\}$.  Recall the function $g(j)$ defined in the proof
of Lemma \ref{lemma2}.  We will only count the paths $\xi\in A_k$
for which $g(\xi_0)=g(\xi_1)=\dots=g(\xi_{j-1})=0$, $g(\xi_j)=i$,
$g(\xi_{j+1})=g(\xi_{j+2})=\dots=g(\xi_{j+m})=d$, and
$g(\xi_{j+m+1})=g(\xi_{j+m+2}=\dots=g(\xi_{n-1})=0$, where $0\leq
j\leq n-(m+1)$ (see Figure \ref{gjpaths}).  The range of $\xi_{n-1}$
is $(n,k_n(\gamma))$. These paths form a subset of the paths which
are counted when computing $\dim(n,k_n(\gamma))$. For a fixed $j$,
the number of such paths is
$$a_0^ja_ia_d^{m}a_0^{n-m-j-1}=a_0^{n-m-1}a_ia_d^m.$$

  \begin{figure}
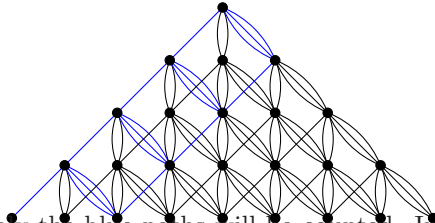

  \begin{center}
\begin{graph}(4,2.5)(0,0)
\unitlength=.7\unitlength
\roundnode{9}(0,0)\roundnode{10}(1,0)\roundnode{11}(2,0)\roundnode{12}(3,0)
\roundnode{13}(4,0)\roundnode{14}(5,0)\roundnode{15}(6,0)\roundnode{4}(1,1)\roundnode{5}(2,1)
\roundnode{6}(3,1)\roundnode{7}(4,1)\roundnode{8}(5,1)\roundnode{1}(2,2)\roundnode{2}(3,2)\roundnode{3}(4,2)
\roundnode{0}(3,3)\roundnode{16}(-1,-1)\roundnode{17}(0,-1)\roundnode{18}(1,-1)\roundnode{19}(2,-1)
\roundnode{20}(3,-1)\roundnode{21}(4,-1)\roundnode{22}(5,-1)\roundnode{23}(6,-1)\roundnode{24}(7,-1)
 \edge{0}{1}[\graphlinecolour(0,0,1)]
\bow{0}{2}{-.1}\bow{0}{2}{.1}
\edge{0}{3}[\graphlinecolour(0,0,1)]\bow{0}{3}{-.1}[\graphlinecolour(0,0,1)]\bow{0}{3}{.1}[\graphlinecolour(0,0,1)]
\edge{1}{4}[\graphlinecolour(0,0,1)]
\bow{1}{5}{-.1} \bow{1}{5}{.1}
\edge{1}{6}[\graphlinecolour(0,0,1)]\bow{1}{6}{-.1}[\graphlinecolour(0,0,1)]\bow{1}{6}{.1}[\graphlinecolour(0,0,1)]
\edge{2}{5}
\bow{2}{6}{-.1} \bow{2}{6}{.1}
\edge{2}{7}\bow{2}{7}{-.1}\bow{2}{7}{.1}
\edge{3}{6}[\graphlinecolour(0,0,1)] \bow{3}{7}{-.1}\bow{3}{7}{.1}
\edge{3}{8}\bow{3}{8}{-.1}\bow{3}{8}{.1}
\edge{4}{9}[\graphlinecolour(0,0,1)]
\bow{4}{10}{-.1} \bow{4}{10}{.1}
\edge{4}{11}[\graphlinecolour(0,0,1)]\bow{4}{11}{-.1}[\graphlinecolour(0,0,1)]\bow{4}{11}{.1}[\graphlinecolour(0,0,1)]
\edge{5}{10}
\bow{5}{11}{-.1} \bow{5}{11}{.1}
\edge{5}{12}\bow{5}{12}{-.1}\bow{5}{12}{.1}
\edge{6}{11}[\graphlinecolour(0,0,1)]
\bow{6}{12}{-.1}\bow{6}{12}{.1} \edge{6}{13} \bow{6}{13}{-.1}
\bow{6}{13}{.1}  \edge{7}{12} \bow{7}{13}{-.1} \bow{7}{13}{.1}
\edge{7}{14} \bow{7}{14}{-.1} \bow{7}{14}{.1} \edge{8}{13}
\bow{8}{14}{-.1} \bow{8}{14}{.1} \edge{8}{15} \bow{8}{15}{-.1}
\bow{8}{15}{.1} \edge{9}{16}[\graphlinecolour(0,0,1)]
\bow{9}{17}{-.1} \bow{9}{17}{.1}
\edge{9}{18}[\graphlinecolour(0,0,1)]\bow{9}{18}{-.1}[\graphlinecolour(0,0,1)]\bow{9}{18}{.1}[\graphlinecolour(0,0,1)]
\edge{10}{17}
\bow{10}{18}{-.1} \bow{10}{18}{.1}
\edge{10}{19}\bow{10}{19}{-.1}\bow{10}{19}{.1}
\edge{11}{18}[\graphlinecolour(0,0,1)]
\bow{11}{19}{-.1}\bow{11}{19}{.1}
\edge{11}{20}\bow{11}{20}{-.1}\bow{11}{20}{.1} \edge{12}{19}
\bow{12}{20}{-.1}\bow{12}{20}{.1}
\edge{12}{21}\bow{12}{21}{-.1}\bow{12}{21}{.1} \edge{13}{20}
\bow{13}{21}{-.1}\bow{13}{21}{.1} \edge{13}{22} \bow{13}{22}{-.1}
\bow{13}{22}{.1}  \edge{14}{21} \bow{14}{22}{-.1} \bow{14}{22}{.1}
\edge{14}{23} \bow{14}{23}{-.1} \bow{14}{23}{.1} \edge{15}{22}
\bow{15}{23}{-.1} \bow{15}{23}{.1} \edge{15}{24} \bow{15}{24}{-.1}
\bow{15}{24}{.1}
\end{graph}
\end{center}
\caption{Only the blue paths will be counted.
In this case $i=d=2$ and $m=0$.}\label{gjpaths}
\end{figure}

Letting $j$ range over $0$ to $n-(m+1)$ we see that
$$\dim(n,k_n(\gamma))\geq (n-m)a_0^{n-m-1}a_ia_d^m.\text{  Hence}$$
$$\begin{array}{rl}
\mu(C)&=\ds\lim_{n\to\infty}\frac{\dim(C,(n,k_n(\gamma)))}{\dim(n,k_n(\gamma))}\\\\
&\leq\ds
\lim_{n\to\infty}\frac{a_0^{n-l}}{(n-m)a_0^{n-m-1}a_ia_d^m}=\lim_{n\to\infty}\frac{a_0^{m+1-l}}{(n-m)a_ia_d^m}\\
&=0.\end{array}$$

By the invariance of $\T$, no cylinder whose terminal vertex is
$(j,k)$, where $j\geq l$, has positive measure.  Using this same
argument for vertices $(j,i)$ where $i\leq k$, we see that the only
edges on which $\mu$ is supported are the edges on the far left of
the diagram.  In this case we have an odometer, and the measure is
as stated in the proposition.  A symmetric argument shows that the
only measure supported on paths for which $dn-k_n(\gamma)\leq k$ for
all $n\geq 0$ is supported on the paths for which $k_n(\gamma)=dn$
for all $n\geq 0$.  Hence $\mu$ is as stated in the proposition.
\end{proof}

\section{Eigenvalues and Total Ergodicity}\label{4.5}

In this section we discuss the eigenvalues of various systems in $\Sl$.
We will show that every Bratteli-Vershik system
determined by a positive integer polynomial of degree 1 for which
either of the coefficients is greater than 1 has at least one non-trivial root of unity as an eigenvalue.  In contrast we will show that the Euler adic has no root of unity (other than one) as an eigenvalue.  In \cite{FP} it is shown that the reverse Euler adic has every root of unity as an eigenvalue.

\begin{theo}  Let $(\X,\T)$ be the Bratteli-Vershik system in $\Slp$ determined
by $p(x)=a_0+a_1x$, with fully-supported, $\T$-invariant, ergodic
probability measure $\mu$.  Then $e^{2\pi i/(a_0a_1)}$, $e^{2\pi
i/a_0}$, and $e^{2\pi i/a_1}$ are eigenvalues of $\T$.\label{eigthm}
\end{theo}

\begin{proof}The sets $\{\gamma\in \X|k_n(\gamma)=0 \text{ for all
}n=0,1,\dots\}$ and $\{\gamma\in \X|k_n(\gamma)=n \text{ for all
}n=0,1,\dots\}$ are both $\T$-invariant sets.  Since $\mu$ is
ergodic and has full support these are sets of measure 0. Recall
that the minimal cylinder
into vertex $(n,k)$,  is denoted $Y_n(k,0)$ and that
$Y_n(k,i)=T^i(Y_n(k,0))$ for $i=0,1,\dots, \dim(n,k)-1$. For
$\mu$-almost every $\gamma\in \X$ there exist $n\geq 0$, $0<k<n$,
and $0\leq j \leq \dim(n,k)-1$ for which $\gamma\in Y_n(k,j)$.

Define the function $f:\X\to\C$ by the following: for $n>1$,
$0<k<n$, and $0\leq j<\dim(n,k)$,
$$f(Y_n(k,j))=(e^{2\pi i/(a_0a_1)})^{j+1},$$
and $f=0$ elsewhere.

 In order to show that $f$ is well defined it is enough to show
that for a minimal cylinder $C$ and an extension $C^j$ of $C$
$$f(C^j)=e^{2\pi i/(a_0a_1)}.$$
We will divide the argument
into two cases,  the case when $C$ is extended by an edge $j\in
\{0,1,\dots,a_1-1\}$ ($C$ extended to the right) and the case when
$C$ is extended by an edge $j\in \{a_1,a_1+1\dots,a_1+a_0-1\}$ ($C$
extended to the left). First assume that $C$ terminates at vertex
$0<k<n$, and let $C^j$ be the extension of the cylinder $C$ by the
edge $j\in \{0,1,\dots, a_1-1\}$, which has terminal vertex
$(n+1,k+1)$. Then $C^j=Y_{n+1}(k+1,j\dim(n,k)).$ Since $0<k<n$,
\begin{equation} \dim(n,k)={n\choose k}a_0^{n-k}a_1^k=0\ \text{mod }a_0a_1,
\label{bicoeffmod}\end{equation} and therefore $$f(C^j)=(e^{2\pi i
/(a_0a_1)})^{j\dim(n,k)+1}=(e^{2\pi i
/(a_0a_1)})^{j\dim(n,k)}e^{2\pi i/(a_0a_1)}=e^{2\pi i/(a_0a_1)}.$$
Now let $C^j$ be the extension of the cylinder $C$ by an edge $j\in
\{a_1, a_1+1,\dots, a_1+a_0-1\}$.  Then
$C^j=Y_{n+1}(k,a_1\dim(n,k-1)+(j-a_1)\dim(n,k))$.  Since $0<k<n$,
$$a_1\dim(n,k-1)={n\choose k-1}a_0^{n-k+1}a_1^{k}=0\ \text{mod }a_0a_1,$$
and from Equation \ref{bicoeffmod} we see that,
\begin{align*}f(C^j)&=(e^{2\pi
i/(a_0a_1)})^{a_1\dim(n,k-1)+(j-a_1)\dim(n,k)+1}\\&=(e^{2\pi
i/(a_0a_1)})^{a_1\dim(n,k-1)}(e^{2\pi
i/(a_0a_1)})^{(j-a_1)\dim(n,k)}e^{2\pi i/(a_0a_1)}\\&=e^{2\pi
i/(a_0a_1)}.\end{align*} Hence $f$ is well defined $\mu$-almost
everywhere.

For $n>1$, $0<k<n$, $0\leq j<\dim(n,k)-1$ and $\gamma\in Y_n(k,j)$,
it is clear that $f(\T \gamma)=e^{2\pi i/(a_0a_1)}f(\gamma)$.  For
$n\geq 0$, $0<k<n$, and $\gamma\in Y_n(k,\dim(n,k)-1)$, there
are $m\geq 0$ and $0<l<m$ such that $\T \gamma\in Y_m(k,0)$.  Then
$f(\gamma)=1$ and $f(\T \gamma)=e^{2\pi i/(a_0a_1)}$.  Hence for
$\mu$-almost every $\gamma\in X$, $f(\T \gamma)=e^{2\pi
i/(a_0a_1)}f(\gamma)$, and $e^{2\pi i/(a_0a_1)}$ is an eigenvalue of
$\T$.

The same argument can be repeated using the eigenvalues $e^{2\pi
i/a_0}$ and $e^{2\pi i/a_1}$.
\end{proof}

\begin{cor}
Let $(\X,\T)$ be the Bratteli-Vershik system determined by the
polynomial $a_0+a_1x$ with a fully-supported, $\T$-invariant,
ergodic probability measure $\mu$.  If either $a_0$ or $a_1$ is
greater than 1, then $\T$ is not weakly mixing.\label{weakmixing}
\end{cor}

\begin{rem} The main result of this section is possible because for
a degree one polynomial, $p(x)=a_0+a_1x$, all the coefficients of
$(p(x))^n$ except the coefficients of $x^0$ and $x^{dn}$ are
divisible by $a_0a_1$. For a polynomial $p(x)$ of degree higher than 1, the
coefficients do not necessarily have a common
factor; therefore this argument is
not sufficient for polynomials of higher degree.
\end{rem}

We will now show that the Euler adic is totally ergodic, in other words,
it has no roots of unity (other than
1) as eigenvalues.

\begin{lemma}  Let $(X,T)$ be the Bratteli-Vershik system determined
by the Euler graph.  Let $\gamma\in X$ be a path such that there is
an $n\in \N$ for which $n-k_n(\gamma)+1>1$,
$k_{n+1}(\gamma)=k_n(\gamma)+1$, and $\gamma_n$ is not the largest
edge (with respect to the edge ordering) connecting
$(n,k_n(\gamma))$ and $(n+1,k_n(\gamma)+1)$. Then
$d(T^{A(n,k_n(\gamma))}\gamma,\gamma)=2^{-n}$.\label{Lemma5.4.1}\end{lemma}

\begin{proof}  For the remainder of the proof let $k=k_n(\gamma)$.
Recall that $Y_n(k,0)$ and $Y_n(k,A(n,k)-1)$ are respectively the
minimal and maximal cylinders into vertex $(n,k)$. There are $l,m\in
\N$ such that $T^l\gamma\in Y_n(k,A(n,k)-1)$, $\gamma\in
T^m(Y_n(k,0))$, and $l+m+1=A(n,k)$.

\begin{figure}[h]
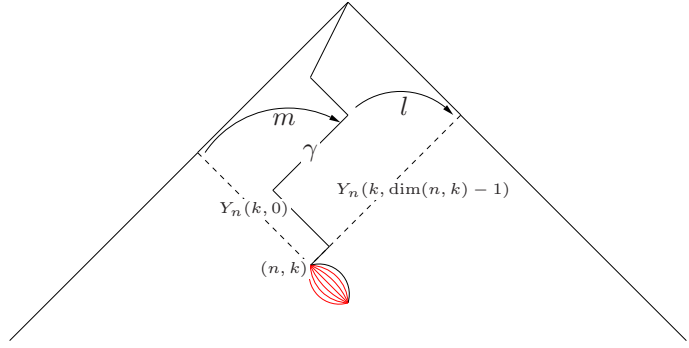

\begin{center}
\begin{graph}(9,4.5)
\unitlength=.5\unitlength \graphnodesize{0}
\roundnode{0}(9,9)\roundnode{1}(5,5)\roundnode{2}(12,6)
\roundnode{3}(8,7)\roundnode{4}(9,6)\roundnode{5}(7,4)\roundnode{6}(8.5,2.5)
\roundnode{7}(8,2)\roundnode{8}(0,0)\roundnode{9}(18,0) \edge{0}{1}
\edge{1}{8}\edge{0}{2}\edge{2}{9}\edge{0}{3}\edge{3}{4}\edge{4}{5}\edge{5}{6}\edge{6}{7}
\edge{1}{7}[\graphlinedash{2}]\edgetext{1}{7}{\tiny $Y_n(k,0)$}
\edge{2}{7}[\graphlinedash{2}]\freetext(11,4){\tiny
$Y_n(k,\dim(n,k)-1)$}
\roundnode{4a}(9.2,6.2)\roundnode{4b}(8.8,5.8)\roundnode{2b}(11.8,6)
\roundnode{1b}(5.2,5) \dirbow{4a}{2b}{.2}\dirbow{1b}{4b}{.2}
\freetext(10.5,6.25){$l$}\freetext(7.3,5.9){$m$} \roundnode{10}(9,1)
\bow{7}{10}{.1}[\graphlinecolour(1,0,0)]
\bow{7}{10}{.2}[\graphlinecolour(1,0,0)] \bow{7}{10}{.3}
\edge{7}{10}[\graphlinecolour(1,0,0)]
\bow{7}{10}{-.1}[\graphlinecolour(1,0,0)]
\bow{7}{10}{-.2}[\graphlinecolour(1,0,0)]
\bow{7}{10}{-.3}[\graphlinecolour(1,0,0)] \edgetext{4}{5}{$\gamma$}
\nodetext{7}(-.7,-.1){\tiny$(n,k)$}
\end{graph}
\end{center}
\caption{$\gamma_n$
is one of the red paths, and $m+l+1=\dim(n,k_n(\gamma))$.}
\end{figure}

$T(T^l\gamma)\in Y_n(k,0)$ and $(T^{l+1}\gamma)_n$ is
the successor of $\gamma_n$ with respect to the edge ordering. Then
$T^m(T^{l+1}\gamma)=T^{A(n,k)}\gamma\in T^m(Y_n(k),0)$. Hence
$$(T^{A(n,k)}\gamma)_0=\gamma_0,\
(T^{A(n,k)}\gamma)_1=\gamma_1,\dots,(T^{A(n,k)}\gamma)_{n-1}=\gamma_{n-1}$$
and $(T^{A(n,k)}\gamma)_n\neq\gamma_n$; therefore
$d(T^{A(n,k)}\gamma,\gamma)=2^{-n}$.
\end{proof}

Key elements of the proof of the following Lemma are already in
\cite{Chacon}, and a similar result has been known for a long time
for substitution and related systems, see \cite{Chacon, Host, Q,
Solo, PS, PetMel, Mela, BDM}.

\begin{lemma}  Let $(X,T)$ be the Bratteli-Vershik system determined
by the Euler graph with the symmetric measure $\eta$. If $\lambda$
is an eigenvalue for $T$, then $\lambda^{A(n,k_n(\gamma))}\to 1$
$\eta$-almost everywhere.\label{convergetoone}
\end{lemma}

There is a closed-form formula for $A(n,k)$ which can be found in
\cite{Salama}:
$$A(n,k)=\sum_{j=0}^{k}(-1)^j {n+2\choose j}(k+1-j)^{n+1}.$$

We will apply the following theorem of Lucas to the above equation.

\begin{theo}[E. Lucas \cite{Lucas}] Let $p$ be a prime number and
$j\leq n$.  Consider the base $p$ decompositions of $n$ and
j:$$n=n_0+n_1p+\dots+n_sp^s$$
$$j=j_0+j_1p+\dots + j_sp^s$$
where $0\leq j_i,n_i<p$ for all i.  Then
$${n\choose j}\equiv_p{n_0\choose j_0}\dots{n_s\choose j_s},$$
with the convention that ${n_i\choose j_i}=0$ if
$j_i>n_i$.\label{lucas}
\end{theo}

\begin{theo}  The Bratteli-Vershik system determined
by the Euler graph with the symmetric measure $\eta$ is totally
ergodic.\label{EulerTotallyErgodic}
\end{theo}
\begin{proof}
It is enough to show that for any prime $p$, $e^{2\pi i/p}$ is not
an eigenvalue for $T$.

If $\lambda=e^{2\pi i/p}$ is an eigenvalue of $T$, by Lemma
\ref{convergetoone} we know that $\lambda^{A(n,k_n(\gamma))}\to 1$
for $\eta$-almost every $\gamma\in X$.  Since $\lambda$ is a root of
unity, for $\eta$-almost every $\gamma$ in $X$, there must be an $N$
such that $n\geq N$ implies $\lambda^{A(n,k_n(\gamma))}=1$.
Therefore for $\eta$-almost every $\gamma\in X$ there is an $N$ such
that $n\geq N$ implies
$$A(n,k_n(\gamma))=0 \mod p.$$

We will show that for every $\gamma\in X$, there are infinitely many
$n$ for which $A(n,k_n(\gamma))\equiv_p1$.  In particular, for every
$l=0,1,\dots$, and $0\leq k\leq p^l-1$,  $A(p^l-1,k)\equiv_p1$.
Recall that for $k\geq 1$,
$$A(p^l-1,k)=\sum_{j=0}^{k}(-1)^j{p^l+1\choose j}(k+1-j)^{p^l}.$$

We will examine this sum by computing each term in the sum
$\text{mod }p$.

\noindent For $j=0$ we have \begin{equation}{p^l+1\choose
0}(k+1)^{p^l}=(k+1)^{p^l}\equiv_p k+1\text{ by Fermat's Little
Theorem}.\label{1Euler}\end{equation} For $j=1$, we have
\begin{equation}(-1){p^l+1\choose 1}k^{p^l}=-(p^l+1)k^{p^l}\equiv_p -k.\label{2Euler}\end{equation}
For $2\leq j\leq p^l-1$ we have $$(-1)^j{p^l+1\choose
j}(k-j)^{p^l}.$$ By Theorem \ref{lucas}, $$(-1)^j{p^l+1\choose
j}(k-j)^{p^l}\equiv_p(-1)^j{1\choose j_0}{0\choose
j_1}\dots{0\choose j_{l-1}}{1\choose 0}(k-j)^{p^l}.$$ Since $2\leq
j\leq p^l-1$, at least one of $j_1,j_2,\dots,j_{l-1}$ must be
positive.  Therefore \begin{equation}(-1)^j{p^l+1\choose
j}(k-j)^{p^l}\equiv_p 0.\label{3Euler}\end{equation}

For fixed $p$ and $0\leq k\leq p^l-1$, we will now compute
$A(p^l-1,k)$.
\begin{align*}
A(p^l-1,0)=(-1)^0{p^l+1\choose 0}(1)^{p^l}=1.\end{align*} Combining
Equations \ref{1Euler}, \ref{2Euler}, and \ref{3Euler} we see that
for $k>0$,
\begin{align*}
A(p^l-1,k)=\sum_{j=0}^{k}(-1)^j{p^l+1\choose
j}(k+1-j)^{p^l}\equiv_p(k+1)-k\equiv_p1.\end{align*}

Hence $\lambda^{A(n,k_n(\gamma))}$ does not converge to 1, and
therefore $T$ has no roots of unity as eigenvalues.
\end{proof}

The status of weak mixing for both the Euler adic and the Pascal adic is still open.

\section{Subshifts and Entropy}\label{3point3}

Recall that for systems in $\Sl$ for which $d=1$, we have that for
all $n=0,1,\dots$ and $k=0,1,\dots n$, $|V_n|=n+1$ and there are
edges between vertices $(n,k)$ and $(n+1,k)$ as well as $(n,k)$ and
$(n+1,k+1)$.

  Let $(X,T)\in \Sl$ with $d=1$. Denote the edges leaving $v_0=(0,0)$ by
  $e_1,e_2,\dots,e_m$.  Define $P_i=\{\gamma\in X|\gamma_0=e_i\}$.  Then
  $\p=\{P_1,P_2,\dots P_m\}$ is a finite partition of $X$ into pairwise disjoint nonempty clopen cylinder sets.
There is a function on $X$, also denoted by $\p$ such that by $\p(\gamma)=j$ for all
$\gamma\in P_j$, $j=1,\dots,m$.  For each $n=0,1,2,\dots$, the
\emph{$\p$-$n$-name} of $\gamma$ is the finite block
$$ \p_0^n(\gamma)=\p(\gamma)\p(T\gamma)\dots\p(T^n\gamma), $$ and the
  \emph{$\p$-name} of $\gamma$ is the doubly infinite sequence
  $$\p_{-\infty}^{\infty}(\gamma)=\dots\p(T^{-2}\gamma)\p(T^{-1}\gamma).\p(\gamma)\p(T\gamma)\p(T^2\gamma)\dots.$$

\begin{figure}[h]
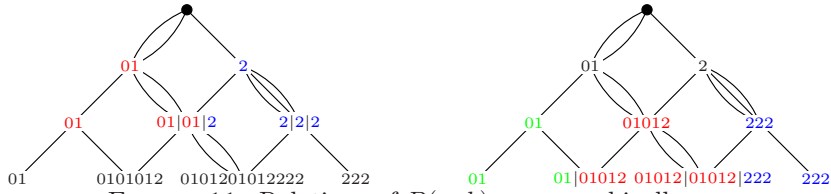

  \begin{center}
      \begin{graph}(6,2)(-1,.5)
      \unitlength=.75\unitlength
  \roundnode{1}(3,3)\roundnode{2}(2,2)\roundnode{3}(4,2)\roundnode{4}(1,1)
  \roundnode{5}(3,1)\roundnode{6}(5,1)\roundnode{7}(0,0)\roundnode{8}(2,0)
  \roundnode{9}(4,0)\roundnode{10}(6,0)\bow{1}{2}{-.1}\bow{1}{2}{.1}\edge{1}{3}
  \edge{2}{4}\bow{2}{5}{-.1}\bow{2}{5}{.1}\edge{3}{5}\edge{3}{6}\bow{3}{6}{-.1}
  \bow{3}{6}{.1}\edge{4}{7}\edge{4}{8}\edge{5}{8}\bow{5}{9}{-.1}\bow{5}{9}{.1}
  \edge{6}{9}\edge{6}{10}
  \nodetext{2}{\tiny\textcolor{red}{01}}\nodetext{3}{\tiny\textcolor{blue}{2}}\nodetext{4}{\tiny\textcolor{red}{01}}
  \nodetext{5}{\tiny$\textcolor{red}{01}|\textcolor{red}{01}|\textcolor{blue}{2}$}
  \nodetext{6}{\tiny$\textcolor{blue}{2}|\textcolor{blue}{2}|\textcolor{blue}{2}$}
  \nodetext{7}{\tiny 01}\nodetext{8}{\tiny 0101012}
  \nodetext{9}{\tiny 0101201012222}\nodetext{10}{\tiny 222}
  \end{graph}
  \begin{graph}(6,2)(-1,.5)
  \unitlength=.75\unitlength
  \roundnode{1}(3,3)\roundnode{2}(2,2)\roundnode{3}(4,2)\roundnode{4}(1,1)
  \roundnode{5}(3,1)\roundnode{6}(5,1)\roundnode{7}(0,0)\roundnode{8}(2,0)
  \roundnode{9}(4,0)\roundnode{10}(6,0)\bow{1}{2}{-.1}\bow{1}{2}{.1}\edge{1}{3}
  \edge{2}{4}\bow{2}{5}{-.1}\bow{2}{5}{.1}\edge{3}{5}\edge{3}{6}\bow{3}{6}{-.1}
  \bow{3}{6}{.1}\edge{4}{7}\edge{4}{8}\edge{5}{8}\bow{5}{9}{-.1}\bow{5}{9}{.1}
  \edge{6}{9}\edge{6}{10}
  \nodetext{2}{\tiny 01}\nodetext{3}{\tiny 2}\nodetext{4}{\textcolor{green}{\tiny 01}}
  \nodetext{5}{\textcolor{red}{\tiny 01012}}\nodetext{6}{\textcolor{blue}{\tiny 222}}
  \nodetext{7}{\textcolor{green}{\tiny 01}}\nodetext{8}{\tiny$\textcolor{green}{01}|\textcolor{red}{01012}$}
  \nodetext{9}{\tiny$\textcolor{red}{\tiny 01012}|\textcolor{red}{\tiny 01012}|\textcolor{blue}{222}$}
  \nodetext{10}{\textcolor{blue}{\tiny 222}}
  \end{graph}
  \end{center}
    \caption{Relations of $B(n,k)$ seen graphically.}\label{concat}
  \end{figure}

  For every vertex
  $(n,k)\in \mathcal{V}$ and $\gamma\in Y_n(k,0)$ define
  $$B(n,k)=\p(\gamma)\p(T\gamma)\p(T^2\gamma)\dots\p(T^{\dim(n,k)-1}\gamma).$$
$B(n,k)$ is called the \emph{basic block at vertex $(n,k)$}.
  Let $l(n,k)$ denote the number of edges connecting $(n,k)$ and
  $(n+1,k)$, and $r(n,k)$ denote the number of edges connecting
  $(n,k)$ and $(n+1,k+1)$.  Then
  \begin{equation}
  B(n+1,k+1)=B(n,k)^{r(n,k)}B(n,k+1)^{l(n,k+1)},
  \label{recursionequation}
  \end{equation}
  where the exponents indicate concatenation, see Figure \ref{concat}.

\begin{defs}
  Let $\Sigma$
  denote the space of bi-infinite sequences on $\{1,2,\dots,m\}$ for which every
  finite subsequence appears as a subblock in some $B(n,k)$,  and let $\sigma : \Sigma\to \Sigma$ denote the
  shift map.\label{definitionofthepartition}
\end{defs}

We define $X'$ to consist of the maximal set $X_{\max}$, its orbit,
and the set of paths that never leave the far left or far right
sides of the diagram: \begin{lemma} Define $X'\subset X$ to consist
of the following paths:\\
1.  $\mathcal{O}(X_{\max});$\\
2. $\{\gamma\in X|k_n(\gamma)=0\  \forall n\in \N\};$.\\
3.  $\{\gamma\in X|k_n(\gamma)=dn\  \forall n\in \N\}$. \\
Then for any fully supported, $T$-invariant, ergodic measure $\mu$,
$\mu(X')=0$.\label{X'}
\end{lemma}
\begin{proof}
  Since $X_{\max}$ is countable, the sets of paths that never leave
  the far left $(k_n\equiv 0)$ or far right $(k_n\equiv n)$ sides of
  the diagram are proper closed $T$-invariant sets, and $\mu $ is
  ergodic, $\mu(X')=0$.
  \end{proof}

\begin{theo}  Let $(X,T) \in \Sl$ with $d=1$, and let $\mu$ be a fully-supported
$T$-invariant ergodic probability measure on $X$.  Let $\Sigma$ be
the subshift defined above.  Then there are a set $X'\subset X$ with
$\mu(X')=0$ and a one-to-one Borel measurable map $\phi: X\setminus
X' \to \Sigma$ such that $\phi\circ T=\sigma\circ \phi$ on
$X\setminus X'$. \label{coding}
\end{theo}

\begin{proof}
For each $\gamma\in X$ define $\phi(\gamma)$ to be the $\p$-name of
$\gamma$. Then for all $\gamma\in X$,
\begin{align*}
\phi\circ T(\gamma) &=\dots \p(T^{-1}\gamma)\p(\gamma).\p(T\gamma)\p(T^2\gamma)\dots\\
&=\sigma (\dots \p(T^{-2}\gamma)\p(T^{-1}\gamma).\p(\gamma)\p(T\gamma)\dots)\\
&=\sigma\circ \phi(\gamma).
\end{align*}

It is clear that $\phi^{-1}$ of any cylinder in $\Sigma$ is a union
of cylinder sets in $X$, hence $\phi$ is Borel measurable.  Defining
$X'$ as above, Lemma \ref{X'} tells us that $\mu(X')=0$.

The strategy for showing $\phi$ is one-to-one is to show that for
$\gamma, \xi\in X\setminus X'$, $\gamma\neq \xi$, there is a
coordinate $j$ such that either $\phi(\gamma)_j$ or $\phi(\xi)_j$ is
a symbol from $B(1,0)$ and the other is a symbol from $B(1,1)$.  This is done in the straightforward but tedious manner of considering
cases according to the different ways that $\gamma,\xi\in
X\setminus X'$ disagree.  We leave the details to the reader.  \end{proof}

\begin{cor}  If $(X,T)$ is a Bratteli-Vershik system in $\Sl$ such that $d=1$,
with fully supported, $T$-invariant, ergodic measure $\mu$ and Borel
sets $\mathcal{B}$, the partition $\mathcal{P}$ is a generating
partition.\label{corollarygeneratingpartition}
\end{cor}

While we only know that this partition is generating on systems in
$\Sl$ with $d=1$, we can define such a partition by the first edge
for any system in $\Sl$, and define $(\Sigma, \sigma)$ as in
Definition \ref{definitionofthepartition}.

\begin{lemma} For large $n$ the number of words of length $n$ appearing in $\Sigma$
is bounded above by a polynomial in $n$ (and hence has topological
entropy 0).\label{1entropy0}
\end{lemma}

\begin{proof}
  Recall that for each vertex $(n,k)$ in $V_n$ and a path
  $\gamma\in Y_n(k,0)$,
  $B(n,k)=\p(\gamma)\p(T\gamma)\dots\p(T^{\dim(n,k)-1}\gamma)$.

  At each level $l$, we determine the maximum possible number of new words
  of length $n$ formed by concatenating two words $B(l,k_1)$ and
  $B(l,k_2)$.  The concatenation of $B(l,k_1)$ and $B(l,k_2)$ can form at
  most $n-1$ new words.  Since there are $dl+1$ vertices, there are
  $(dl+1)^2$ possible distinct concatenations.  Hence there are at
  most $(dl+1)^2(n-1)$ new words formed by concatenation.

At level $n$ all blocks except possibly $B(n,0)$ and $B(n,dn+1)$
have length at least $n$.  Concatenating $B(n,0)$ and $B(n,1)$
creates the word $B(1,0)^n$.  For all levels $m\geq 0$, the edges of
the diagrams dictate that $B(m,0)$ only joins with $B(m,1)$, hence
the concatenation of $B(m,0)$ and $B(m,1)$ at levels $m\geq n$ will
only create $B(1,0)$ across their juncture and hence no longer
create words that have not been seen before. Likewise for
$B(n,dn-1)$,  $B(n,dn)$, and $B(1,d)^n$. All other blocks at level
$n$ are of length at least $n$. Since all words in subsequent levels
are created by some concatenations of entries on level $n$, no more
new words are formed.

Therefore the number of words of length $n$ is bounded above by
\begin{align*}\sum_{l=1}^n(dl+1)^2(n-1)
&\leq n^2(dn+1)^2\\
&\leq d^2n^4+2dn^3+n^2.
\end{align*}
\end{proof}

\begin{theo}  Let $(X,T)$ be a Bratteli-Vershik system in $\Sl$ with
a $T$-invariant measure $\mu$.  Then $(X,T,\mu)$ has entropy
0.\end{theo}

\begin{proof}
One may replace the partition in the proof of Theorem
\ref{1entropy0} by the partition by the first $l$ edges, $\p_l$, and
use the subshift $(\Sigma_l,\sigma)$ corresponding to the partition
$\p_l$. A similar counting argument will yield that the number of
$n$-blocks in $\Sigma_l$ is again bounded by a polynomial in $n$.
Now for each $n=1,2,\dots,$
\begin{align*}H_{\mu}\left(\bigvee_{i=0}^{n-1}T^i\p_l\right)&=-\sum_{A\in
\bigvee_{i=0}^{n-1}T^i\p_l}\mu(A)\log(\mu(A))\\
& \leq H_{ud}\left(\bigvee_{i=0}^{n-1}T^i\p_l\right),
\end{align*}
where $ud$ is the measure on $\bigvee_{i=0}^{n-1}T^i\p_l$ that gives
each element equal measure.  If $p_n$  is the cardinality of
$\bigvee_{i=0}^{n-1}T^i\p_l$, then
\begin{align*}
H_{ud}\left(\bigvee_{i=0}^{n-1}T^i\p_l\right)&=-\sum_{i=1}^{p_n}\frac{1}{p_n}\log\left(\frac{1}{p_n}\right)\\
&=\log(p_n).
\end{align*}
Since there is a constant $c_l$ such that $p_n\leq n^{c_l}$ for all
large $n$, the entropy of the system with respect to the partition
$\p_l$ is
$$h_\mu(\p_l,T)=\lim_{n\to
\infty}\frac{1}{n}H_\mu\left(\bigvee_{i=0}^{n-1}T^i\p_l\right)\leq
\lim_{n\to\infty}\frac{c_l}{n}\log(n)=0.$$

 Now let $\B_l$ be
the $\sigma$-algebra generated by $\p_l$ and $\B$ the
$\sigma$-algebra of $(X,T)$.  Since $\B_l\nearrow \B$, $h_\mu(T)$ is
the limit of $h_\mu(\p_l,T)$. Hence the entropy of $(X,T,\mu)$ is 0.
\end{proof}

\section{Loosely Bernoulli}\label{looselybernoullipolynomial}
 The property of loosely Bernoulli was introduced by Feldman in
 \cite{Feldman} as well as by Katok and Sataev in \cite{KS}.
 A transformation that has zero entropy (see \cite{PetBook}) is loosely
 Bernoulli if and only if it is isomorphic to an induced map of an irrational
 rotation on the circle.

\begin{defs}  The $\overline{f}$ distance between two words $v=v_1\dots v_l$ and $w=w_1\dots w_l$
of the same length $l>0$ on the same alphabet is
$$\overline{f}(v,w)=\frac{l-s}{l},$$
where $s$ is the greatest integer in $\{0,1,\dots,l\}$ such that
there are $1\leq i_1<i_2<\dots <i_s\leq l$ and $1\leq j_1<j_2<\dots
<j_s\leq l$ with $v_{i_r}=w_{j_r}$ for $r=1,\dots,s$.
\end{defs}

\begin{defs} Let $T$ be a zero-entropy measure-preserving transformation on the probability space
$(X,\mathcal{B},\mu)$, and let $\mathcal{P}$ be a finite measurable
partition of $X$.  Then the process $(\mathcal{P},T)$ is said to be
\emph{loosely Bernoulli} (LB) if for all $\varepsilon>0$ and for all
sufficiently large $l$ we can find $A\subset X$ with
$\mu(A)>1-\varepsilon$ such that for all $\gamma,\xi\in A$,
$$
\overline{f}(\mathcal{P}_0^l(\gamma),\mathcal{P}_0^l(\xi))<\varepsilon.$$
$T$ is said to be \emph{loosely Bernoulli} if $(\p,T)$ is loosely
Bernoulli for all partitions $\p$.
\end{defs}

$T$ is LB if for a generating partition $\p$, $(\p,T)$ is LB.  Some
of the Bratteli-Vershik systems determined by positive integer
polynomials have already been shown to be loosely Bernoulli.
Janvresse and de la Rue proved it for the Pascal adic in \cite{JdlR},
and in \cite{Mela}, M\'ela showed it for polynomials of arbitrary
degree where all the coefficients are 1. We have established this
property for Bratteli-Vershik systems determined by arbitrary
positive integer polynomials as well as for the Euler adic.

\begin{theo} The Bratteli-Vershik systems $(\X,\T)$ in $\Slp$ determined by positive integer polynomials
are loosely Bernoulli with respect to each of their $\T$-invariant
ergodic probability measures.\label{polylooseber}
\end{theo}

The proof of Theorem \ref{polylooseber} will follow the ideas of
Janvresse and de la Rue.  There are two cases, depending on whether
or not the ergodic measure has full support. The following lemma
gives a seemingly weaker sufficient condition for the loosely
Bernoulli property to hold.

\begin{lemma}[Janvresse, de la Rue \cite{JdlR}] Let $(X,\mathcal{B},\mu,T)$ be a measure-preserving
system with entropy 0.  Suppose that for every $\varepsilon
>0$ and for $\mu\times\mu$-almost every $(\gamma,\xi)\in X\times X$ we
can find an integer $l(\gamma,\xi)\geq 1$ such that
$$\overline{f}\left(\mathcal{P}_0^{l(\gamma,\xi)}(\gamma),\mathcal{P}_0^{l(\gamma,\xi)}(\xi)\right)<\varepsilon.$$
Then the process $(\p,T)$ is LB.
\end{lemma}

\begin{lemma}  Let $(\X,\T)$ be the Bratteli-Vershik system determined by the positive integer polynomial
$p(x)$ of degree $d$, with ergodic, $\T$-invariant probability
measure $\mu$. For $\mu\times\mu$-almost every $(\gamma,\xi)$ in
$\X\times \X$, we can find arbitrarily large $n$ such that
$k_n(\gamma)=k_n(\xi)$. \label{meet}
\end{lemma}

\begin{proof}
Define the random variables $\{Z_i\}_{i=1}^n$ from $\X$ to
$\{0,1,\dots,d\}$ by letting
$Z_i(\gamma)=k_i(\gamma)-k_{i-1}(\gamma)$. This is an i.i.d.
process. Define $S_n:\X\times\X\to \{-d,-d+1,\dots,0,1,\dots,d\}$ by
$S_n(\gamma,\xi)=\sum_{i=1}^n[Z_i(\gamma)-Z_i(\xi)]$.  $(S_n)$ is a
symmetric random walk, and hence recurrent.  Thus for
$\mu\times\mu$-almost every $(\gamma,\xi)$ in $\X\times\X$, there
are infinitely many $n$ such that $k_n(\gamma)=k_n(\xi)$.
\end{proof}

\begin{prop}  Let $(\X,\T)$ be a Bratteli-Vershik system determined
by a positive integer polynomial of degree $d$, with a
fully-supported, $\T$-invariant, ergodic probability measure $\mu$.
Let $\p$ be the partition determined by the first edge.  Then the
process $(\p,\T)$ is loosely Bernoulli. \label{LBprop}
\end{prop}
\begin{proof}
The partition $\mathcal{P}$ was described in detail in Section
\ref{3point3} for polynomials of degree 1.  The same notations will
be used here.  In particular, recall that for each vertex $(n,k)$,
and a path $\gamma\in Y_n(k,0)$,
$$B(n,k)=\p(\gamma)\p(\T \gamma)\dots\p(\T^{\dim(n,k)-1}\gamma).$$

From Lemma \ref{meet}, we know that for $\mu\times\mu$-almost every
$(\gamma,\xi)\in \X\times\X$, $\gamma$ and $\xi$ meet infinitely
often. Also, for $\mu$-almost every $\gamma$, for each cylinder $C$,
$$\frac{\dim(C,(n,k_n(\gamma)))}{\dim(n,k_n(\gamma))}\to \mu(C).$$
Let $p_0=\mu([0])$ denote the weight associated to each edge
labeled 0. For each $r\geq 1$, with probability $p_0^{2r}>0$, both
$\gamma$ and $\xi$ continue along edges labeled $0$ for the next $r$
edges.

Given $\varepsilon>0$, choose $r$ so that $p_0^r<\varepsilon/2$. Let
$C$ be a cylinder with terminal vertex $(r,dr)$.  Then for
$(\mu\times\mu)$-almost every $(\gamma,\xi)$ there are infinitely
many $n$ for which $k_n(\gamma)=k_n(\xi)=k$, $\ds \left|
\frac{\dim(C,(n,k_n(\gamma)))}{\dim(n,k_n(\gamma))}-\mu(C)\right|<\frac{\varepsilon}{2}\text{,
and}
\gamma_n=\xi_n=\gamma_{n+1}=\xi_{n+1}=\dots=\gamma_{n+r-1}=\xi_{n+r-1}=0$.
Then the $\p$-names of both $\gamma$ and $\xi$ have long central
block $B(n+r,k+dr)$.  If we decompose $B(n+r,k+dr)$ into blocks from
level $n$, we see that the first block to appear is $B(n,k)$.  Both
$\gamma$ and $\xi$ have their decimal point in this first block of
$B(n,k)$.

\begin{figure}[h]
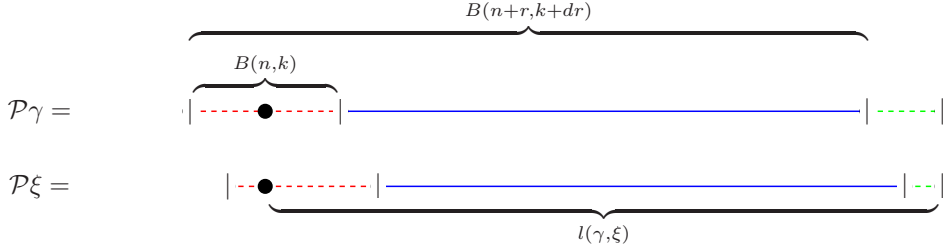

\begin{center}
\begin{graph}(12,3)(-2,.5)
\roundnode{1}(0,2)\roundnode{2}(1,2)\roundnode{3}(2,2)\roundnode{4}(9,2)\roundnode{5}(10,2)
\roundnode{6}(.5,1)\roundnode{7}(1,1)\roundnode{8}(2.5,1)\roundnode{9}(9.5,1)\roundnode{10}(10,1)
\nodetext{1}{$|$}\edge{1}{2}[\graphlinedash{2}\graphlinecolour(1,0,0)]
\edge{2}{3}[\graphlinedash{2}\graphlinecolour(1,0,0)]\edge{3}{4}[\graphlinecolour(0,0,1)]
\edge{4}{5}[\graphlinecolour(0,1,0)\graphlinedash{2}]\edge{6}{7}[\graphlinecolour(1,0,0)\graphlinedash{2}]
\edge{7}{8}[\graphlinecolour(1,0,0)\graphlinedash{2}]\edge{8}{9}[\graphlinecolour(0,0,1)]
\edge{9}{10}[\graphlinedash{2}\graphlinecolour(0,1,0)]
\nodetext{3}{$|$}\nodetext{4}{$|$}\nodetext{5}{$|$}\nodetext{6}{$|$}\nodetext{8}{$|$}\nodetext{9}{$|$}
\nodetext{10}{$|$}
\freetext(-2,2){$\p \gamma=$}\freetext(-2,1){$\p \xi=$}
\freetext(4.5,3.2){$\overbrace{\hskip 3.55in}^{B(n+r,k+dr)}$}
\freetext(1,2.5){$\overbrace{\hskip .75in}^{B(n,k)}$} \freetext(5.5,
0.5){$\underbrace{\hskip 3.5in}_{l(\gamma,\xi)}$}
\end{graph}
\end{center}
\caption[$\p \gamma$ and $\p \xi$ agree on a long block.]{$\p
\gamma$ and $\p \xi$ agree on the blue line, which is the rest of
$B(n+r,k+dr)$ after the end of the initial $B(n,k)$.}
\end{figure}

If we let $l(\gamma,\xi)=|B(n+r,k+dr)|$ we have
$$\overline{f}(\p_0^{l(\gamma,\xi)}\gamma,\p_0^{l(\gamma,\xi)}\xi)\leq
\frac{|B(n,k)|}{|B(n+r,k+dr)|}=\frac{\dim(n,k)}{\dim(n+r,k+dr)}.$$
By the isotropic nature of the diagram,
$\dim(n,k)=\dim(C,(n+r,k+dr))$. Hence,
$$\frac{|B(n,k)|}{|B(n+r,k+dr)|}=\frac{\dim(n,k)}{\dim(n+r,k+dr)}=\frac{\dim(C,(n+r,k+dr))}{\dim(n+r,k+dr)}.$$
By Lemma \ref{lemma1},
$$\frac{\dim(C,(n+r,k+dr))}{\dim(n+r,k+dr)}\to \mu(C)=p_0^r.$$
We can now take $n$ large enough so that
$$\frac{|B(n,k)|}{|B(n+r,k+dr)|}<p_0^r+\frac{\varepsilon}{2}<\varepsilon.$$
Hence $(\p,\T)$ is loosely Bernoulli.
\end{proof}

If the partition by the first edge is a generating partition, then
Proposition \ref{LBprop} would be enough to say that all
Bratteli-Vershik systems in $\Slp$ are loosely Bernoulli.  As
it stands, we are able to prove that all the systems in $\Slp$ are
loosely Bernoulli without proving that the partition by the first edge is generating.

\begin{cor}  Let $(\X,\T)$ be a Bratteli-Vershik system determined
by a positive integer polynomial $p(x)$ of degree $d$ with
fully-supported, $\T$-invariant, ergodic probability measure $\mu$.
Let $\p_l$ be the partition determined by the first $l$ edges.  Then
the process $(\p_l,\T)$ is loosely Bernoulli. \label{LBcor}
\end{cor}
\begin{proof}
By telescoping to the levels which are multiples of $l$, the
Bratteli-Vershik system becomes the system determined by the
polynomial $q(x)=(p(x))^l$, and $\p_l$ becomes the partition on the
first edge. By Proposition \ref{LBprop}, this system is LB.
\end{proof}

\begin{theo}[Ornstein, Rudolph, Weiss \cite{ORW}]  If $G$ is a compact group and $\phi:G\to G$
is rotation by $\rho$ with $\Z_\rho$ dense in $G$, then for any
partition $\p$, $(\p,\phi)$ is LB. \label{ORW1} \end{theo}

\begin{theo}[Onrstein, Rudolph, Weiss \cite{ORW}]  If $\mathcal{B}_n\nearrow \mathcal{B}$ and
$(X,\mathcal{B}_n,T)$ is LB for each $n$, so is
$(X,\mathcal{B},T)$.\label{ORW2}
\end{theo}

\begin{proof}[Proof of Theorem \ref{polylooseber}]
  If $\mu$ does not have full support, by Proposition \ref{nfsthm},
  the Bratteli-Vershik system is a stationary odometer.  Every odometer is a
  compact group rotation, and hence by Theorem \ref{ORW1} is LB.

  For a fully-supported measure $\mu$, Corollary \ref{LBcor} says
  that for each $l$, the process $(\p_l,\T)$ is LB.  Let
  $\mathcal{B}_l$ be the $\sigma$-algebra generated by $\p_l$.  Then
  $(\X,\mathcal{B}_l,\mu,\T)$ is LB and $\mathcal{B}_l\nearrow
  \mathcal{B}$.  Hence Theorem \ref{ORW2} tells us that
  $(\X,\mathcal{B},\mu,\T)$ is LB.
  \end{proof}

We now give the same result for the Euler adic.

\begin{theo}  Let $(X,T)$ be the Bratteli-Vershik system in $\Sl$
determined by the Euler graph and with the symmetric measure $\eta$.
Then $T$ is loosely Bernoulli.\label{EulerAdicLB}\end{theo}

\begin{proof}
  Let $\p$ be the partition according to the first edge, described in Section \ref{3point3} for Bratteli-Vershik
  systems in $\Sl$ for which $d=1$.  By Corollary
  \ref{corollarygeneratingpartition} this is a generating partition;
  therefore it is sufficient to show that the process $(\p,T)$ is LB.

Proposition 2 in \cite{BKPS} tells us that for
$\eta\times\eta$-almost every $(\gamma,\xi)\in X\times X$,
$k_n(\gamma)=k_n(\xi)=k$ infinitely many times.  For such an $n$,
with conditional probability
$$\left(\frac{n-k+1}{2n+2}\right)^2,$$
$k_{n+1}(\gamma)=k_n(\gamma)+1=k_{n+1}(\xi)$ and both $\gamma_n$ and
$\xi_n$ are one of the first $(n-k+1)/2$ edges into $(n+1,k+1)$.
Lemma 2 in \cite{BKPS} tells us that $k_n(\gamma)/n\to 1/2$ $\eta$-almost everywhere.
Therefore, for $\eta$-almost every $\gamma\in X$,
$$\frac{n-k_n(\gamma)+1}{2n+2}\to\frac{1}{4}.$$
Then for $\eta$-almost every $\gamma\in X$ we can take $n$ large
enough so that
$$\left(\frac{n-k_n(\gamma)+1}{2n+2}\right)^2>\frac{1}{64}.$$ Hence the set
of $(\gamma,\xi)$ for which there are infinitely many $n$ such that
$k_n(\gamma)=k_n(\xi)=k$, $k_{n+1}(\gamma)=k_{n+1}(\xi)=k+1$, and
each of $\gamma_n$ and $\xi_n$ are one of the first $(n-k+1)/2$
edges connecting $(n,k)$ and $(n+1,k+1)$ has full measure.

Then the $\p$-names of both $\gamma$ and $\xi$ have long central
block $B(n+1,k+1)=B(n,k)^{n-k+1}B(n,k+1)^{k+2}$.  Both $\gamma$ and
$\xi$ have their decimal point in this first block of $B(n,k)$. For
some subblocks $w_0w_1\dots w_{j_1}$ and $w_0'w_1'\dots w_{j_2}'$ of
$B(n,k)$ and $m_1,m_2$ with $(n-k+1)/2\leq m_1,m_2\leq n-k+1$,
\[\begin{array}{l}\p_0^{\infty}\gamma=\dots_\bullet w_0w_1\dots
w_{j_1}(B(n,k))^{m_1}(B(n,k+1))^{k+2}\dots\\
\p_0^{\infty}\xi=\dots_\bullet w_0'w_1'\dots
w_{j_2}(B(n,k))^{m_2}(B(n,k+1))^{k+2}\dots\end{array}\]

Then $\p_0^{\infty}\gamma$ and $\p_0^{\infty}\xi$ agree on
$\min\{m_1, m_2\}$ consecutive blocks $B(n,k)$.  Let
$l(\gamma,\xi)=\min\{m_1,m_2\}A(n,k)+\max\{j_1,j_2\}$. Then
$$\overline{f}(\mathcal{P}_0^{l(\gamma,\xi)}(\gamma),\mathcal{P}_0^{l(\gamma,\xi)}(\xi))\leq
\frac{\max\{j_1,j_2\}}{\min\{m_1,m_2\}|B(n,k)|}=
\frac{2A(n,k)}{(n-k+1)A(n,k)}=\frac{2}{n-k+1}.$$

Lemma 2 in \cite{BKPS} says that $k_n(\gamma)/n\to 1/2$ as
$n\to \infty$.  Thus given $\varepsilon>0$, we can let $n$ be large
enough so that
$$\frac{2}{n-k+1}<
\varepsilon.$$  Then $T$ is LB.
\end{proof}

\bibliographystyle{plain}
\bibliography{bibliography}
\end{document}